\pdfoutput=1
\RequirePackage{ifpdf}
\ifpdf 
\documentclass[pdftex]{sigma}
\else
\documentclass{sigma}
\fi

\usepackage{tcolorbox}
\tcbuselibrary{breakable}

\usepackage{pgf,tikz,pgfplots,tikz-3dplot}
\usetikzlibrary{calendar,folding}
\usetikzlibrary{arrows,patterns,decorations.pathmorphing,backgrounds,positioning,fit,petri}
\usetikzlibrary{calc,3d,intersections,shapes}
\usepgfplotslibrary{polar}
\usetikzlibrary {arrows.meta}
\usetikzlibrary{shapes.misc}

\numberwithin{equation}{section}

\newtheorem*{Theorem*}{Theorem}

{ \theoremstyle{definition}

}

\begin{document}

\allowdisplaybreaks

\renewcommand{\thefootnote}{}

\newcommand{\arXivNumber}{2307.16234}

\renewcommand{\PaperNumber}{066}

\FirstPageHeading

\ShortArticleName{Fragments of a History of the Concept of Ideal}

\ArticleName{Fragments of a History of the Concept of Ideal.\\ Poncelet's and Chasles's Reflections on Generality\\ in Geometry and their Impact on Kummer's Work\\ with Ideal Divisors\footnote{This paper is a~contribution to the Special Issue on Differential Geometry Inspired by Mathematical Physics in honor of Jean-Pierre Bourguignon for his 75th birthday. The~full collection is available at \href{https://www.emis.de/journals/SIGMA/Bourguignon.html}{https://www.emis.de/journals/SIGMA/Bourguignon.html}}}

\Author{Karine CHEMLA}

\AuthorNameForHeading{K.~Chemla}

\Address{Laboratoire SPHERE UMR 7219, Universit\'e Paris Cit\'e,\\
 27~rue Jean-Antoine de Ba\"{\i}f, 75013 Paris, France}
\Email{\href{mailto:chemla@univ-paris-diderot.fr}{chemla@univ-paris-diderot.fr}}
\URLaddress{\url{http://www.sphere.univ-paris-diderot.fr/spip.php?article78}}

\ArticleDates{Received June 30, 2023, in final form June 27, 2024; Published online July 21, 2024}

\Abstract{In this essay, I argue for the following theses. First, Kummer's concept of ``ideal prime factors of a complex number'' was inspired by Poncelet's introduction of ideal elements in geometry as well as by the reconceptualization that Michel Chasles put forward for them in 1837. In other words, the idea of ideal divisors in Kummer's ``theory of complex numbers'' derives from the introduction of ideal elements in the new geometry. This is where the term ``ideal'' comes from. Second, the introduction of ideal elements into geometry and the subsequent reconceptualization of what was in play with these elements were linked to philosophical reflections on generality that practitioners of geometry in France developed in the first half of the 19${}^{\rm th}$ century in order to devise a new approach to geometry, which would eventually become projective geometry. These philosophical reflections circulated as such and played a key part in the advance of other domains, including in Kummer's major innovation in the context of number theory.}


\Classification{01A55; 00A30}

\begin{flushright}
\begin{minipage}{50mm}
\it For Jean-Pierre Bourguignon,\\
 as a token of friendship
 \end{minipage}
 \end{flushright}

\renewcommand{\thefootnote}{\arabic{footnote}}
\setcounter{footnote}{0}

\section{Introduction}\label{section1}

The first piece of evidence we have about Ernst Eduard Kummer's (1810--1893) introduction of ``ideal prime factors of a complex number'' (\textit{ideale Primfactoren einer complexen Zahl}) in number theory is a letter he wrote to Leopold Kronecker (1823--1891) on October 18, 1845.\footnote{Kurt Hensel (1861--1941) selected parts of the letters that Kummer sent to Kronecker between 1842 and 1865 for publication in \cite[p.~46--102]{45}. By ``parts'', I mean in particular that Hensel deleted whole passages of the letters, replacing them by the sign ``\dots''. The letter from October 18, 1845, is on pp.~64--68. Unfortunately, Hensel did not reproduce it entirely. The original document seems to be lost \cite{19}. The whole volume in which these letters were published was reprinted in Andr\'e Weil's edition of Kummer's works \cite[p.~30--133]{46}. The letter in question is reprinted in \cite[p.~94--97]{46}. See also \cite[p.~5]{46}. In what follows, I will quote those of Kummer's publications that were reprinted in this context on the basis of this edition.} The letter clearly shows that Kummer has just made a breakthrough and was writing to his correspondent in order to put his ideas in order \cite[p.~94]{46}. Four months later, in March 1846, the theory, in the context of which ideal factors were introduced, had taken shape and was made public for the first time, through an ``excerpt from his latest researches in number theory'' that Johann Peter Gustav Lejeune-Dirichlet (1805--1859) communicated to the Berlin Academy at Kummer's request \cite[p.~87]{23}. The following year, this publication was reprinted with some changes, in Crelle's \textit{Journal}, under the title ``Zur Theorie der complexen Zahlen''.\footnote{\cite{25}, reprinted in \cite[p.~203--210]{46}. The meaning of ``complex numbers'' in this title will be clarified in the next section. In \cite{25}, Kummer's first communication to the Berlin Academy was erroneously said to date from 1845. The mistake was pointed out by Edwards in \cite[p.~224]{17} and by Weil in \cite[p.~4]{46}. The changes that occurred between the 1846 and 1847 versions have not yet been systematically described. In this article, I will point out some of them when they occur in the passages of the article that I quote. An English translation of the version of the article that appeared in Crelle's \emph{Journal} was published decades ago \cite{12}. Edwards \cite[p.~224, fn~20]{17} pointed out that it is ``flawed''. In the present article, I quote substantial parts of \cite{23}~and I provide a new English translation for them. I plan to publish a complete English translation of this text, which will indicate systematically the differences between the 1846 and the 1847 versions.} Two features of this first publication (\cite{23}, reprinted in 1847 as \cite{25}) are essential for the purpose of this article.\looseness=1

 Firstly, in March 1846, Kummer published only an overview of his theory, without any proofs. This might be due to the fact that this article reflects an academic presentation, in which Kummer simply gave a non-technical presentation of his main ideas. Except for issues of chronology, historiography has paid less attention to this text than to later, more complete, publications that Kummer went on to make on the same topic. However, I will argue that this overview sheds some light on the process by which Kummer was able to achieve his breakthrough. Several kinds of clues support this view -- notably, as we shall see, a coherent set of terms that Kummer used as well as the structure of his article. In contrast to his March 1846 communication, some months later, in \cite{24} -- dated September 1846 -- Kummer presented his theory in greater detail, and, some years later, in 1851 \cite{26}, he published it with the intention of giving ``an entire and continuous treatise'' (``Trait\'e entier et continu'').\footnote{\cite[p.~378]{26}, reprinted in \cite[p.~363--484]{46}. The quotation comes from p.~364.~Note that \cite{24}, reprinted in \cite[p.~211--251]{46}, was published in Crelle's \emph{Journal} immediately after the above-mentioned article \cite{25}.} These two publications do not reflect as clearly as the March 1846 presentation the elements that Kummer drew his inspiration from for his introduction of ideal divisors -- even though, as will be shown in Section~\ref{section5.3}, the 1851 memoir still contains key clues.\looseness=1

 Secondly, to situate what he meant by ``ideal complex numbers'' in the context of the mathematics of his time, in the March 1846 publication, Kummer sketched several comparisons. He presented the uniform decomposition into prime factors that his ideal factors allowed him to achieve as being underlain by the same idea as the uniform decomposition of polynomials in a single variable into linear factors, which was made possible by the introduction of complex numbers (in the usual current sense of the term, which, as we shall shortly see, differed from Kummer's use of the expression). Kummer also compared the need for his new type of factors with the need that had led Gauss to introduce complex numbers of the type $a+bi$ ($a$~and~$b$ being integers) in his research about biquadratic residues. These comparisons are not surprising, and historical studies have mainly dwelled on them. Historiography has also given pride of place to the comparison between the factorization of a ``complex number'' into prime (sometimes ideal) factors and the decomposition of chemical substances into elements. In fact, Kummer introduced the latter comparison only later (our earliest evidence for this comparison dates from June 1846). He did so to illustrate the procedure by which he made ideal factors visible and how they combined with one another into actual complex numbers.\footnote{See the subsequent extant letter that Kummer wrote to Kronecker on June 14, 1846 \cite[p.~98]{46}. We return to this letter below. The analogy with chemistry recurs in \cite[p.~359--361]{24} and in \cite[p.~447--448]{26}, see also, respectively, \cite[p.~243--245, 433--434]{46}. } However, the first publication of March 1846 puts forward a third comparison, which, to my knowledge, has remained unnoticed until our collective work on generality \cite{4, 9}.\footnote{I should mention a passing remark in \cite[p.~165]{1}. However, Avigad does not analyze what lies behind this remark. Moreover, he only associates the reference to geometry with Poncelet, whereas I will argue that the reference to Chasles's understanding of ideal elements in geometry was also important for Kummer.} Indeed, shifting the point of the comparison to the issue of the \textit{definition} of his new entities, Kummer drew a comparison between his ``true (usually ideal) prime factors'' and the ``chord common to two circles'', which is ``ideal'' when the circles do not intersect -- referring to the search of an ``actual definition of the ideal common chord'' in the latter circumstances (``eine wirkliche Definition dieser idealen gemeinschaftlichen Sehne'') \cite[p.~88]{23}.

 The expression ``ideal common chord'' unmistakably evokes Jean-Victor Poncelet (1788--1867), who introduced it, notably in his \textit{Trait\'e des propri\'et\'es projectives des figures} (1822)~-- a~book that played a major role in the emergence of projective geometry. Given that Kummer, in his first publication on the subject, makes this comparison explicit, it is striking that he chose to refer to his new type of ``prime factors of a complex number'' by precisely the same term ``ideal'' that Poncelet had used some two decades earlier. This choice of terminology, as well as Kummer's explicit parallel in 1846 between his approach and the new geometry that French mathematicians like Gaspard Monge (1746--1818), Jean-Victor Poncelet and Michel Chasles (1793--1880) had contributed to develop, both point to the same conclusion. Both seem to suggest that this new geometry played a part in the considerations that led Kummer to introduce his ``ideal complex numbers''.

In the present article, I argue that the new geometry played a role in the thought process that Kummer followed in his introduction of ideal divisors, and, also, that this role is most obvious in the March 1846 presentation. The reason for this is not only that in the latter article, Kummer makes explicit the parallel between the two types of ideal elements in relation to the issue of forming the correct definition (feature 2 of the 1846 presentation that I described above), but also that the structure of the March 1846 communication is intimately related to the way the text emphasizes the issue of definition (feature 1 of the 1846 presentation that was underlined). Indeed, the article does this in a way that also stems from the new geometry. In other words, the emphasis placed on definition is correlated with Kummer's reference to the ``ideal'' elements introduced by the new geometry. To my knowledge, the reference to this geometry disappears from Kummer's later publications on the topic, except for the 1851 publication in French (\cite[p.~430]{26}, see also \cite[p.~416]{46}) in which it recurs. Interestingly, in the context of the re-occurrence, Kummer again insisted on the issue of defining ideal elements. As I argue below (Section~\ref{section5.3}), this latter text yields evidence converging towards the same conclusion about the impact that the new geometry had on Kummer's conception of ideal divisors, even though the 1851 article no longer mentions the use of the term ``ideal'' in geometry.\looseness=1

 In brief, in this essay, I argue for the following theses. First, Kummer's concept of ``ideal prime factors of a complex number'' was inspired both by Poncelet's introduction of ideal elements in geometry and by Michel Chasles's reconceptualization of them in 1837. That is, the idea of ideality in Kummer's ``theory of complex numbers'' has its origins in the introduction of ideal elements in the new geometry. This is where the term ``ideal'' comes from. Second, the introduction of ideal elements in geometry and the subsequent reconceptualization of what was in play with these elements were linked to philosophical reflections on generality that practitioners of geometry in France developed in order to devise a new approach to geometry -- which would eventually become projective geometry. These philosophical reflections circulated as such and played a key part in the advance of other domains, including in Kummer's major innovation in the context of number theory. I claim that Kummer was aware of Poncelet's work on ideal elements -- I return to this issue below -- and, also, that he read Chasles's writings on the topic, which were instrumental for the development of his new theory.

More precisely, my argument will rely in an essential way on the following facts: The structure of the 1846 publication introducing these ideal factors, as well as the set of terms about definitions Kummer used in it in relation to ideal elements, yield striking evidence about how Kummer progressively shaped his definition of ``ideal prime factors''. What is more, the elements of this process that we can restore highlight the fact that Kummer followed Chasles's recommendations on how to elaborate general and uniformly valid definitions in the new geometry, instead of introducing Poncelet's ``ideal elements''. In other words, the 1846 publication shows that, to provide the uniform generality in the decomposition of complex numbers into prime factors that Kummer sought to achieve, he followed a strategy inspired by actors who had been involved in the development of the new geometry. This is important because, after all, the key point of Kummer's breakthrough lies in the definition he created for ideal elements. Establishing these facts will enable us to see that the comparison that Kummer drew between ideal elements in both contexts is an essential one, which explains how and why the origin of the term ``ideal'' lies in the new geometry. This case study thus illustrates how mathematicians' historical and philosophical reflections can be crucial for their mathematical work.

 To establish these theses, I will mainly rely on the four pieces of Kummer's writings that I~have mentioned above: his two letters to Kronecker from October 18, 1845, and from June 14, 1846; the March 1846 text communicated to the Berlin Academy; and the 1851 publication in French.

\section{The first public introduction of ``ideal complex numbers''}\label{section2}

 Let us begin by examining in greater detail how, in 1846, Kummer \cite[p.~87]{23} introduced his results and his motivations, in the first public and formal presentation of his ``theory of complex numbers'' -- an expression that also formed the title of the publication of this text in 1847~\cite{25}.

Kummer defined the notion of a ``complex number'' that is central to his work as early as in the second paragraph of his presentation \cite[p.~87--88]{23}. For him, if $\alpha$ denotes an imaginary root of the equation $\alpha^\lambda = 1$, in which $\lambda$ is a prime integer, a ``complex number'' $f(\alpha) $ -- or alternatively a ``complex integral number (complexe ganze Zahl)'' -- $ $is an expression of the following kind:\looseness=1
\[
f(\alpha)= a + a_1\alpha + a_2\alpha ^2+\dots +a_{\lambda -1}\alpha ^{\lambda -1},
\]
in which the $a_{i}$ are integers. Gauss had already done some work on these numbers, which, in modern terms, are called cyclotomic integers.

Kummer's project, the completion of which is announced in the 1846 presentation, aims to decompose any complex number of this type (and notably prime integers) into prime factors, and this in a unique way. In other words, Kummer aimed at developing a theory of divisibility for these numbers identical to that holding for common integers. However, to achieve this goal Kummer identified a key difficulty, which he depicted as follows:\footnote{The quotation follows the version in \cite[p.~88]{23}. My emphasis.}

\renewenvironment{quote}{
 \list{}{%
 \leftmargin1.0cm
 \rightmargin0cm
 \topsep=2pt
 }
 \item\relax
}
{\endlist}

\begin{quote}
Such a complex number can either be decomposed into factors of the same kind, or also cannot.
In the first case, this number is composite, whereas in the other, the number was \textit{until now called a complex prime number}. However, I have now noticed that, even if $f(\alpha)$ \textit{can in no way be decomposed into complex factors}, it does \textit{not yet have the true nature of a complex prime number}, because it already usually \textit{lacks} the \textit{first and most important property of prime numbers}: that is, that the product of two prime numbers is not divisible by any prime number different from them. (Eine solche complexe Zahl kann entweder in Factoren derselben
Art zerlegt werden oder auch nicht, im ersten Falle ist sie eine zusammengesetzte,\footnote{Here, \cite[p.~319]{25} has an additional ``Zahl''. I do not indicate the differences in punctuation between the two publications.} im anderen Falle ist sie \textit{bisher eine complexe Primzahl genannt} worden. Ich habe nun aber bemerkt, dass wenn auch $f(\alpha)$ \textit{auf keine Weise in complexe Factoren zerlegt} \textit{werden kann}, sie darum\footnote{Here, instead of ``darum'', \cite[p.~319]{25} has ``deshalb'' (``by way of consequence'').} \textit{noch nicht die wahre Natur einer complexen Primzahl} hat, weil sie schon \textit{der ersten und wichtigsten Eigenschaft der Primzahlen} gew\"ohnlich\footnote{In \cite[p.~319]{25}, ``gew\"ohnlich'' features after ``schon''.} \textit{ermangelt}, n\"amlich dass das Product zweier Primzahlen durch keine von ihnen verschiedene Primzahl theilbar ist.)
\end{quote}

In other words, unlike prime numbers in ordinary arithmetic, which serves here as a model for the development of the theory, for complex numbers, Kummer identified the need to distinguish between indecomposability and primality -- a distinction that, he asserted, had been so far overlooked. It is precisely to bridge the gap between these two notions in the case of complex numbers -- and thereby achieve a uniform decomposition -- that Kummer suggested introducing ``ideal complex numbers''. Indeed, with respect to these indecomposable complex prime numbers that nevertheless do not have the property of prime numbers pointed out, Kummer goes on:

\begin{quote}
\textit{Even though} numbers $f(\alpha)$ of this type \textit{cannot be decomposed} into complex factors, they have nevertheless \textit{much more} the \textit{nature} of \textit{composite numbers}; but the \textit{factors} are consequently \textit{not} \textit{actual}, but \textit{ideal complex numbers}. (Es haben vielmehr solche Zahlen $f(\alpha)$, \textit{wenngleich} sie \textit{nicht} in complexe Factoren \textit{zerlegbar} sind, \textit{dennoch} die \textit{Natur der zusammengesetzten Zahlen}, die \textit{Factoren} aber sind alsdann \textit{nicht wirkliche}, sonder\footnote{Sic. This typo is corrected into ``sondern'' in \cite[p.~319]{25}.} 
\textit{ideale complexe Zahlen}.) (\cite[p.~88]{23}, with my emphasis, except for the last expression, which the author also emphasizes.)
\end{quote}

I translate ``wirkliche {\dots} complexe Zahlen'' here by ``actual {\dots} complex numbers''. By this term, Kummer means complex numbers that can be expressed in the form introduced above $\bigl(f(\alpha)=a+a_1\alpha + a_2\alpha ^2+\dots +a_{\lambda -1}\alpha ^{\lambda -1}\bigr)$. On this basis, one thus understands why Kummer opens his communication to the Berlin academy with the following words:

 I have succeeded in essentially \textit{completing} and \textit{simplifying} the \textit{theory of the complex numbers} that are made up of higher roots of the unity ({\dots}); and, to be precise, I did this by introducing a very specific type of \textit{imaginary divisors}, which I call \textit{ideal complex numbers} ({\dots}) (Es ist mir gelungen \emph{die Theorie derjenigen complexen Zahlen}, welche aus h\"oheren Wurzeln der Einheit gebildet sind, ({\dots}), wesentlich\footnote{The adverb wesentlich (significantly, essentially) is deleted in \cite{25}.} zu \textit{vervollst\"andigen} und zu \textit{vereinfachen}, uud\footnote{Sic, and corrected into ``und'' in \cite{25}.} zwar durch Einf\"uhrung einer ganz\footnote{The adverb ganz (very, wholly) is deleted in \cite{25}.} eigenth\"umlichen Art \emph{imagin\"arer Divisoren}, welche ich \emph{ideale complexe Zahlen} nenne ({\dots}).)\footnote{\cite[p.~87]{23}. Again, my emphasis, except for the expression ``ideal complex numbers''.}

It is noteworthy that, for Kummer, the double effect of introducing ``ideal complex numbers'' is to ``complete'', and to ``simplify'' the theory -- we return to these values and to Kummer's interpretation of them in Section~\ref{section5.3}. More importantly for us, note that these ``ideal complex numbers'' are introduced as ``divisors''. In other words, it is not directly, but indirectly -- through a specific relation that they have with other \textit{actual} numbers -- that they are brought to light. This detail will prove meaningful. In addition, it is interesting to note that Kummer first referred to these ``divisors'' as ``imaginary'', but then discarded this term to introduce the expression of ``ideal complex numbers'' instead. In fact, as we shall see, these ``ideal complex numbers'' are \textit{not} imaginary, if by this term, we mean ``complex numbers'' in the sense in which the expression is used today. What this tells us about Kummer's notion of the imaginary would be interesting to pursue, but I will not follow this line of inquiry here. I will rather focus on what the choice of the term ``ideal'' indicates.
 Indeed, Kummer first justified the introduction of such ``complex numbers'' through two comparisons. He wrote:

\begin{quote}
The same \textit{simple} idea underlies the introduction of such ideal complex numbers and the introduction of imaginary formulas in algebra and analysis, namely in the decomposition of integral rational functions (KC: i.e., polynomials) into their \textit{simplest} factors -- the linear ones. Moreover, it is also the same need through which, necessarily, \emph{Gauss} first introduced the complex numbers of the form $a + b \sqrt{-1}$ in the investigations on the biquadratic remainders, because in this case, all prime factors of the form $4m +1$ displayed the nature of composite numbers.\footnote{The formulation is awkward. However, I have tried to reflect the German language as I understand it here.} (Der Einf\"uhrung solcher idealen complexen Zahlen liegt derselbe \emph{einfache} Gedanke zu Grunde als\footnote{Here, instead of ``als'', Kummer \cite[p.~319]{25} has ``wie \dots''.} der Einf\"uhrung der imagin\"aren Formeln in die Algebra und Analysis, namentlich bei der Zerf\"allung der ganzen rationalen Functionen in die\footnote{Here, instead of ``die'', Kummer \cite{25} has ``ihre''. I translate according to the more recent version.} \emph{einfachsten} Factoren, welche die line\"aren sind.\footnote{Here, instead of ``welche die line\"aren sind'' Kummer \cite{25} has ``die line\"aren''.} Ferner ist es auch dasselbe Bed\"urfniss, durch welches gen\"othigt\footnote{Here, in \cite{25}, Kummer added a comma.} \textit{Gauss} bei den Untersuchungen \"uber die biquadratischen Reste, weil hier alle Primfactoren von der Form $4m +1$ die Natur zusammengesetzter Zahlen zeigten,\footnote{What is here between commas was put between brackets in \cite[p.~319]{25}.} die complexen Zahlen von der\footnote{Here, Kummer \cite[p.~319]{25} added the noun ``Form''. I translate according to the more recent version.} $a + b\sqrt{-1}$ zuerst einf\"uhrte.) \cite[p.~88, my emphasis except for the name of Gauss]{23}.
\end{quote}

 The two examples of the introduction of new entities that Kummer put forward as comparable to his innovation have in common with his ``ideal complex numbers'' the facts that, in the former case, ``the \textit{simplest} factors'' were identified, and, in the latter, numbers that might have been conceived as prime were in fact ``composite''. These two comparisons are followed immediately by a third one, in which precisely the term ``ideal'' is introduced. Let us first read it, before pointing out the key clues that the text contains:
\begin{quote}
 To achieve now a \textit{lasting definition} of the \textit{true} (\textit{usually ideal}) prime factors of the complex numbers, it was \textit{necessary to find out the properties} of the prime factors of complex numbers \textit{that would persist/remain in all circumstances}, which would be absolutely \textit{independent} of the \textit{contingency/accidental circumstances} of whether the \textit{actual} decomposition takes place or not, \textit{more or less precisely as, when in geometry}, one speaks of the \textit{chord common to two circles} \textit{also} when the \textit{circles do not intersect} each other, \textit{one looks for an actual definition of this ideal common chord that fits all the situations of the circles}.

(Um nun zu einer \textit{festen Definition} der \textit{wahren} (\textit{gew\"ohnlich idealen}) Primfactoren der complexen Zahlen zu gelangen, war es n\"othig, \textit{die unter}\footnote{Here, in 1846, Kummer \cite{23} has ``nnter'', which is corrected into ``unter'' in \cite[p.~320]{25}.} \textit{allen Umst\"anden bleibenden Eigenschaften} der Primfactoren complexer Zahlen \textit{aufzusuchen},\footnote{In \cite[p.~320]{25}, this verb ``aufzusuchen'' (to seek) is rewritten into ``zu ermitteln'' (to identify, to determine).} welche von der \textit{Zuf\"alligkeit}, ob die \textit{wirkliche} Zerlegung Statt habe, oder nicht, ganz \textit{unabh\"angig} w\"aren: \textit{ohnegef\"ahr}\footnote{In \cite[p.~320]{25}, this adverb is rewritten into ``ungef\"ahr''.} \textit{ebenso wie wenn in der Geometrie} von der \textit{gemeinschaftlichen Sehne zweier Kreise} ge\-spro\-chen wird, \textit{auch}\footnote{In \cite{25}, ``dann'' is added here, which adds emphasis ``also then \dots''.} wenn \textit{die Kreise sich nicht schneiden,} \textit{eine wirkliche Definition dieser idealen gemeinschaftlichen Sehne gesucht wird,}\footnote{In \cite[p.~320]{25}, the structure of the sentence ``ebenso wie wenn in der Geometrie von der gemeinschaftlichen Sehne zweier Kreise gesprochen wird, (\dots) eine wirkliche Definition dieser idealen gemeinschaftlichen Sehne gesucht wird \dots'' is rewritten into ``eben so, wie man, wenn in der Geometrie von der gemeinschaftlichen Sehne zweier Kreise gesprochen wird, (\dots ) eine wirkliche Definition dieser idealen gemeinschaftlichen Sehne sucht''.} \textit{welche f\"ur alle Lagen der Kreise passt}.) \cite[p.~88, my emphasis]{23}.\looseness=-1
\end{quote}

As I have explained in the introduction, this third comparison that Kummer drew has been largely overlooked in the secondary literature. However, it contains key clues as to the sources on which Kummer relied in order to introduce the ``ideal complex elements'', and also as to the reflections that inspired the procedure he followed to this effect.

 Firstly, Kummer here established a parallel between his ``ideal complex numbers'' and the ``ideal chord common to two circles'' -- a concept introduced by Jean-Victor Poncelet in relation to a reflection on how to increase generality in geometry. The use of the same term ``ideal'' in geometry, which was put forward precisely in the first public presentation of the theory and does not appear in the context of any of the other comparisons mentioned, is certainly striking. Kummer further outlined how he understands the comparison. The situations in which the decomposition of complex numbers into prime factors ``actually'' ``takes place''~-- namely, one can exhibit ``actual'' complex numbers that are the prime factors~-- correspond, in geometry, to the situations in which ``the two circles intersect''~-- namely, one can exhibit ``actual'' points of intersection, which define the common chord. By contrast, the situations in which the decomposition is not ``actual''~-- that is, for which ideal factors need to be introduced~-- correspond, in geometry, to the situations in which two circles do not (actually) intersect -- that is, cases in which, despite the absence of real intersection points, Poncelet still considers that there is a common chord, but an ``ideal'' one. To understand this parallel better, in the next section (Section~\ref{section3}), we will outline the main characteristics of Poncelet's treatment of ideal elements in his \textit{Trait\'e des propri\'et\'es projectives des figures}~\cite{36}.\looseness=1

 Secondly, it is no less striking that in this third comparison, the emphasis is placed on an issue different from what was pointed out in the first two comparisons. Indeed, now, Kummer lay stress on how to define ``the true (usually ideal) prime factors''. For this, he insisted that the definition he looked for should be one that would be uniformly valid for all possible prime factors. Indeed, to achieve this definition, he made explicit that he was looking for ``\textit{the properties} of the prime factors of complex numbers \textit{that would persist/remain in all circumstances}, which would be absolutely \textit{independent} of the \textit{contingency/accidental circumstances} of whether the \textit{actual} decomposition takes place or not''. Note again that, the introduction of ``ideal complex numbers'' is wholly dependent on the part they play, as \textit{factors}, in the theory of divisibility. In Section~\ref{section4}, I will show that this way of shaping definitions for situations such as the one Kummer was facing, as well as the related set of terms Kummer used here, both derive from specific pages of Michel Chasles's 1837 \textit{Aper\c{c}u Historique.} In those pages, Chasles got rid of Poncelet's notion of ideality and, as an alternative, offered a way of dealing with ``ideal common chords'' that is different from Poncelet's.

 As a result, we will see that Kummer's ``theory of complex numbers'' has roots in both Poncelet's and Chasles's reflections on generality in geometry, which he synthesized.

\section{Poncelet's introduction of ``ideal'' elements in geometry}\label{section3}

Let us focus now on the meaning that Poncelet gave to ideal elements, to begin with in his \textit{Trait\'e des propri\'et\'es projectives des figures}, published in 1822 (see Section~\ref{section3.2}). To address this issue, and also to examine how Kummer had access to Poncelet's ideas on these elements, we will consider the treatise in the broader context of Poncelet's geometrical writings.\looseness=1

\subsection{Poncelet's geometrical manuscripts and publications}\label{section3.1}

At the end of his life, between 1862 and 1866, Poncelet published or republished most of his geometrical works in four volumes: two of these volumes were titled \textit{Applications d'analyse et de g\'eom\'etrie qui ont servi de principal fondement au Trait\'e des propri\'et\'es projectives des figures} (volume~1 appeared in 1862 \cite{38} and volume~2 in 1864 \cite{39}), and the other two, \textit{Trait\'e des propri\'et\'es projectives des figures; ouvrage utile \`a ceux qui s'occupent des applications de la g\'eom\'etrie descriptive et d'op\'erations g\'eom\'etriques sur le terrain} (respectively, in 1865~\cite{40} and 1866~\cite{41}). I~will provide some information on these publications here, insofar as this will be useful in determining which of Poncelet's writings were available to Kummer and his contemporaries before 1845. These will also make clear which documents are available to us in order to understand the formation of Poncelet's conception of ideal elements that was instrumental for Kummer.\footnote{For a more detailed biographical sketch (with all necessary references) as well as a first approach to the development of Poncelet's ideas about ideal elements, see \cite{2}. In this other publication, we also give a larger bibliography and an outline of the historiography of ideal elements. I will focus here on what is useful in relation to the topic of the present article.}

The first volume of works that Poncelet published in 1862 \cite{38} consists chiefly of seven notebooks that he had never published before and that he had said he had written in Saratoff, during the year he spent there as prisoner of war. Indeed, on November 18, 1812, while Poncelet was participating in the Russian campaign as a member of the Grande Arm\'ee, he was taken prisoner and held at Saratoff from March 1813 to June 1814. To the manuscripts that Poncelet made public in 1862, he added documents and notes, some of which had been written while he was a~student at the Ecole Polytechnique between 1807 and 1810.\footnote{In \cite{2}, the approach to Poncelet's introduction of ideal elements draws mainly on a comparison between the Saratoff notebooks and the 1822 treatise. Our goal is to shed light on the mathematical problems to which Poncelet constantly returned, notably between 1813 and 1822, and on how his reflections on these problems and his changing solutions are correlated with his introduction of notably ``ideal secants'' and ``ideal chords''. In this other publication, we discuss the distinction between ``secants'' and ``chords'', and its mathematical importance for Poncelet. As I show below, both Chasles and Kummer speak of ``ideal chords'', but in fact, by this term, they refer to ``ideal secants''. I will not dwell on the opposition here.}

The war ended in June 1814. Poncelet could thus leave Saratoff, and he reached Metz in September of the same year. In the second volume of the same publication (which appeared in 1864), Poncelet further made available unpublished drafts, written between 1814 and 1820, in which we see the mathematician's approach to geometry maturing, in what spare time was left to him by his life as an officer in Metz. To these manuscripts, Poncelet added articles that had appeared prior to the publication of his 1822 treatise, as well as documents attesting notably to discussions with fellow mathematicians.

The third volume \cite{40} offers a new edition of the 1822 \textit{Trait\'e des propri\'et\'es projectives des figures,} ``revised, corrected and with new annotations'' (``revue, corrig\'ee et augment\'ee d'annotations nouvelles''), whereas in the fourth volume \cite{41}, Poncelet included annotated versions of articles and memoirs that had been published after the treatise -- notably in the \textit{Journal f\"ur die reine und angewandte Mathematik}, which August Crelle edited in Berlin at the time -- as well as other documents and reflections.

When does the expression ``ideal chord'' mentioned by Kummer first occur in Poncelet's writings? Despite scattered occurrences of this term, as well of the expression of ``ideal secant'' in \cite{38},\footnote{See \cite{2}. We do not know how exactly Poncelet published his notebooks, and in particular how much change he brought to them. Accordingly, we need to be cautious when expressions such as these occur in these notebooks without being used in a systematic way.} as far as I can tell the first work in which Poncelet foregrounds these two notions -- even putting them to the fore in the structure of the text -- is the memoir that he presented to the Academy of Sciences on May 1, 1820 \cite[p.~365]{40}. Its text was published as the fifth notebook in \cite[p.~365--454]{40}, under the title ``Essai sur les propri\'et\'es projectives des sections coniques (Essay on the projective properties of conic sections)'', and the first section of this memoir was entitled precisely ``Notions pr\'eliminaires sur les cordes ou s\'ecantes id\'eales des sections coniques (Preliminary notions on ideal chords and secants of conic sections)''. However, when Kummer presented his own theory in 1846, none of these texts were yet published. Ironically, it seems that Augustin-Louis Cauchy's report on this 1820 memoir, which, to Poncelet's great despair, criticized a key part of his memoir to which we shall return, was the first publication to discuss the notion of an ``ideal chord''~\cite{5}.

Poncelet's 1820 memoir was part of a first draft of the 1822 \textit{Trait\'e des propri\'et\'es projectives des figures} \cite[p.~365]{40}. The treatment of ``ideal chords'' and ``ideal secants'' that the memoir puts forward testifies to ideas about ideality that are quite similar to those presented in the treatise. Let us thus first examine Poncelet's project in this treatise and the role it devotes to ideal elements.

\subsection{Ideal elements in the 1822 \emph{Trait\'e des propri\'et\'es projectives des figures}}\label{section3.2}

As is well known, and as Poncelet explains in the ``Introduction'' to his major work in geometry, his main project grew out of key questions: Why, he asked, is analysis able to proceed uniformly and produce general results, when applied to geometrical issues, and why, by contrast, is ``the ordinary or ancient geometry'' deprived of the same properties? How could one change geometry to make it enjoy these properties?\footnote{\cite[p.~xix--xx]{36}. The third manuscript that was published in \cite[p.~167--295]{39} under the title ``Sur la loi des signes de position en g\'eom\'etrie, la loi et le principe de continuit\'e (On the law of signs of position in geometry and on the law and principle of continuity)'' was written during the winter 1815--1816. This notebook (cahier) nevertheless has marks of changes that date to 1817. In it, Poncelet's reflection seems to bear rather on how to use ``algebraic analysis'' and ``rational geometry'' conjointly (see notably \cite[p.~294]{39}) -- Incidentally, the expression of ``g\'eom\'etrie rationnelle (rational geometry)'' appears to have been used for the first time by Charles Dupin: in the preface to his work, Dupin \cite[p.~viii]{16} gives the expression as synonymous with ``pure geometry'' and opposed to ``analytic geometry''. Nabonnand \cite[p.~18]{33} also notes this point. The article \cite{35}, published on November~1, 1817, appears to me to follow the same line of thought as the third manuscript. To be sure, this was what Monge encouraged students at the Ecole Polytechnique to do \cite{3}. I argue that, by contrast, we can note a shift in Poncelet's reflections with the fourth manuscript (``cahier'') entitled ``Consid\'erations philosophiques et techniques sur le principe de continuit\'e dans les lois g\'eom\'etriques (Philosophical and technical considerations on the principle of continuity in geometrical laws)'' \cite[p.~296--362]{39}. Indeed, in contrast with the previous manuscript, this one, which Poncelet dates from the winter 1818--1819, starts with the statement of a project that is now quite similar to the project motivating the 1822 treatise \cite[p.~296--298]{39}. Friedelmeyer \cite[p.~87--91]{20} appears to consider that~\cite{35} already attests to this project. As we will see, this shift is correlated with a new and massive use of the adjective ``ideal'', which has, however, a meaning different from the meaning that appears in 1820 and will be the one used in the 1822 treatise.} The treatise thus aims at offering general means to ``improve rational geometry'' and notably to increase significantly the generality that can be achieved by proceeding geometrically.\footnote{\cite[p.~xxxiii]{36}. In a text written in 1823, presented at the Acad\'emie des sciences in March 1824, and published in Crelle's \emph{Journal}, \cite{37}, the geometer returns to these principles, putting forward additional principles, and adding further explanations.}

A first general method of this kind, to which Poncelet gave pride of place in the very title of the 1822 treatise, consists in concentrating on the properties of a figure that are preserved by central projection -- these are the properties he calls ``projective''. Indeed, properties of this kind can be proved by reference to a particular figure, and they can subsequently be stated with full generality for any figure that derives from this particular one by such a projection. Poncelet specified that these properties can be ``descriptive'' (that is, depending only on the situations of lines with respect to each other) as well as ``metrical''.\footnote{\cite[p.~xxii, xxiv--xxv]{36}. In \cite{37}, this principle is promoted to become the first one, while in (Poncelet 1822), it features as the second one. In \cite{37}, Poncelet discusses further both the descriptive and metrical properties of this kind, devoting a ``second note'' to the latter. Poncelet's systematic consideration of metrical projective properties has been largely ignored in the secondary literature \cite{2}. Chemla \cite{8}, Friedelmeyer \cite{20} and Nabonnand \cite{32, 33} are exceptions. Lorenat \cite{30} addresses this anachronistic reading of Poncelet's work from a historical perspective.}

In contrast with the latter kind of reasoning, which proceeds from the particular to the general, Poncelet introduced a second principle, which will enable one to conclude that a property established for a general state of a figure holds true for any other state of the same figure.\footnote{On this opposition, and the relationship between Lazare Carnot's principle of correlation between figures and Poncelet's principle of continuity, see \cite{8, 9}. Lorenat's presentation of the principle of continuity in \cite[p.~91]{29} seems to conflate the two principles (see also p.~96, \cite[p.~89]{27} and \cite[p.~176]{28}).} When introducing this principle in his 1822 treatise, Poncelet referred to it as ``the \textit{principle} or \textit{the law of continuity} of mathematical relations of the abstract and figurate magnitude (le \textit{principe} ou \textit{la loi de continuit\'e} des relations math\'ematiques de la grandeur abstraite et figur\'ee)''.\footnote{\cite[p.~xxiii]{36}. Emphasis is his. Nabonnand \cite[p.~3, 17]{33} also emphasizes this part of the name of the principle.} As we will see, the full expression is meaningful, but, to refer to the principle in a more convenient way, we will shorten it below -- as is most often done -- into ``the principle of continuity''. The introduction of ideal elements is intimately related to the reflections Poncelet devoted to this principle between 1815 and 1822. In the third and fourth notebooks, composed between 1815 and 1818 (see footnote 30), Poncelet discussed, under the same name, principles of this kind. Such a principle is also central to the 1820 memoir -- the fifth manuscript in \cite{38} -- and it is precisely the focus of Cauchy's criticism in his report. We leave for a~future publication a~complete examination of the evolution of the principle, with its transformations that can be identified through Poncelet's manuscripts,\footnote{Compare Friedelmeyer in \cite[p.~65--72, 114--116]{20}, Nabonnand in \cite[p.~12--20]{32}, and in \cite[p.~17--22]{33}, which offer some views on this issue.} concentrating here simply on what sheds light on Poncelet's notion of ``ideal elements'' that were to be taken up by Kummer.

In the introduction to his 1822 treatise, Poncelet devoted many pages to a discussion of the ``\textit{principle} or \textit{the law of continuity} of mathematical relations of the abstract and figurate magnitude''. The principle draws in an essential manner on a way of conceiving a figure, with which Poncelet began as follows:

\begin{quote}
Let us consider an \textit{arbitrary} \textit{figure}, in a \textit{general} and somehow \textit{indeterminate} \textit{position} among all those that the figure could take without \textit{violating the laws}, the \textit{conditions} and the \textit{link}s that \textit{subsist} between the various \textit{parts} of the \textit{system} (Consid\'erons une \textit{figure} \textit{quelconque}, dans une \textit{position} \textit{g\'en\'erale} et en quelque sorte \textit{ind\'etermin\'ee}, parmi toutes celles qu'elle peut prendre \textit{sans violer} les \textit{lois}, les \textit{conditions}, la \textit{liaison} qui \textit{subsistent} entre les diverses \textit{parties} du \textit{syst\`eme}) \cite[p.~xxii, my emphasis]{36}.
\end{quote}

Here, Poncelet introduced two terms to refer to the object that is taken as his subject of inquiry. At the beginning of the quotation, he refers to this object as a ``figure'', which can have several ``positions''. The second half of the quotation makes clear what a ``position'' is, by reference to the object now taken as a ``system'', which is composed of parts -- for instance, a~line and a curve in a~plane. Indeed, in the case of this example, the two parts of the system can intersect each other or not: these are two ``general'' positions of the figure, in contrast to a position in which the line would be tangent to the curve. In other words, the figure for Poncelet is not the figure of the ancient geometry, but it is composed of parts and encompasses any configuration that the different parts can take.\footnote{In the 1820 memoir, Poncelet \cite[p.~387]{39} referred to the ``intention that one has to extend the conception of a figure, presently geometric and possible, to all the states that this figure can take, even to those in which certain objects lose their absolute existence (la volont\'e qu'on a d'\'etendre la conception d'une figure, actuellement g\'eom\'etrique et possible, \`a tous les \'etats que peut prendre cette figure, m\^eme \`a ceux o\`u certains objets perdent leur existence absolue.)'' Poncelet (1822: 50) repeated this sentence almost verbatim, with two exceptions: instead of ``all the states that a figure can take'', he writes ``the various states through which a figure can go''. In addition, instead of ``certain objects lose their absolute existence'', he specifies ``their absolute and real existence''.} However, Poncelet added a key restriction, when he specified that, in any position considered, ``the laws, the conditions and the links that subsist between the various parts of the system'' should not be ``violat[ed]''.\footnote{Earlier -- for instance in the fourth manuscript of winter 1818--1819 -- Poncelet \cite[p.~300]{39} had spoken of ``comparing a figure with all those that may be supposed to result from it through the progressive and continuous motion of some of the parts that enter in the figure, without violating the links and the dependence that were primitively established between these parts (de comparer une m\^eme figure avec toutes celles qui peuvent \^etre cens\'ees en r\'esulter par le mouvement progressif et continu de certaines des parties qui y entrent, sans violer la liaison et la d\'ependance primitivement \'etablies entre elles.)'' The same expression recurs in the 1820 memoir (Poncelet 1820: 451).} What did he mean by these terms?

For the example of the curve and the line in a plane precisely, Poncelet illustrated what he meant by ``conditions'' in a letter to Olry Terquem (1782--1862) dated October 14, 1819. In this letter, he replied to objections made to his principle of continuity. Indeed, on November~23, 1818, Poncelet had communicated an early version of his reflections on this principle to three friends with whom he discussed his mathematical ideas, namely, Olry Terquem, Fran\c{c}ois-Joseph Servois (1767--1847) and Charles-Julien Brianchon (1783--1864) \cite[p.~530--539]{39}. In the reply that Terquem sent Poncelet in the name of the three correspondents in September 1819,\footnote{Del Centina \cite[Section~2]{15} dwells on this exchange of letters.} Terquem and Servois put forward the following objection: instead of moving a line with respect to a curve in the plane in which they both lie, one might by continuity suppose that the line is ``lifted out of the plane'' and then ``all the intersections'' between the line and the curve ``become impossible'' \cite[p.~543]{39}. Poncelet's answer points out that when one introduces the possibility that the line moves out of the plane, one changes the ``essential and constitutive conditions defining the system (les conditions essentielles et constitutives du syst\`eme)''. The fact that a line is in a plane requires a ``specific condition (condition distincte)''. ``One violates [this condition] or one changes the question if one then supposes the line to lie in space (on la viole ou {\dots} on change la question, en la supposant ensuite dans l'espace)'' \cite[p.~549]{39}.

Poncelet \cite[p.~95]{36} further gave the example of a ``link'', when he mentioned the relation of involution that exists between the six points that a transversal makes when it passes through the system of a conic section and an inscribed quadrilateral. Note that this link is conceived as persisting, even when the intersections are no longer real.\footnote{This idea is expressed quite clearly, e.g., in the fourth manuscript of the winter 1818--1819 \cite[p.~341--342]{39}.}

It is precisely for figures conceived in this way that the principle of continuity is meaningful. In fact, it even underlies this very conception of the figure. Immediately afterwards, Poncelet formulated the principle as follows:

\begin{quote}
Let us suppose that, on the basis of these data, one has \textit{found one or more relations} or \textit{properties}, either \textit{metrical,} or \textit{descriptive}, that belong to the figure, using the ordinary explicit reasoning, i.e., using the way of proceeding that, in some cases, is considered the only rigorous one. Is it not obvious that if, while preserving the same data, one makes the \textit{primitive} \textit{figure vary by imperceptible degrees}, or one applies to certain parts of this figure a \textit{continuous motion} (in fact arbitrary), is it not obvious that the properties and relations that were \textit{found for the first system} will \textit{continue to apply} to its \textit{successive states}, provided, however, that one pays attention to the particular changes that may have occurred in the system, such as when some magnitudes have vanished, or have changed direction or sign, etc. -- changes that will always be easy to recognize \textit{a priori} and using reliable rules? (supposons que, d'apr\`es ces donn\'ees, on ait \textit{trouv\'e une ou plusieurs relations ou propri\'et\'es}, soit \textit{m\'etriques}, soit \textit{descriptives}, appartenantes (sic) \`a la figure, en s'appuyant sur le raisonnement explicite ordinaire, c'est-\`a-dire par cette marche que, dans certains cas, on regarde comme seule rigoureuse. N'est-il pas \'evident que si, en conservant ces m\^emes donn\'ees, on vient \`a faire varier la \emph{figure primitive par degr\'es insensibles}, ou qu'on imprime \`a certaines parties de cette figure un \emph{mouvement continu} d'ailleurs quelconque, n'est-il pas \'evident que les propri\'et\'es et les relations, \emph{trouv\'ees pour le premier syst\`eme}, \emph{demeureront applicables aux \'etats successifs} de ce syst\`eme, pourvu toutefois qu'on ait \'egard aux modifications particuli\`eres qui auront pu y survenir, comme lorsque certaines grandeurs se seront \'evanouies, auront chang\'e de sens ou de signe, etc., modifications qu'il sera toujours ais\'e de reconna\^{i}tre \textit{\`a priori} et par des r\`egles s\^ures?) \cite[p.~1822: xxii, my emphasis, except for ``metrical'', ``descriptive'' and ``a priori'', which are emphasized by Poncelet]{36}.
\end{quote}

The principle thus asserts that, having established properties or relations through reliance upon a general position of a figure, conceived along the lines discussed above, the same properties or relations ``continue to apply'', even when the figure is in another position. I use the verb ``continue'' to translate here Poncelet's expression, which literally reads: ``remain applicable''.\footnote{Lorenat \cite[p.~96]{29} uses the same expression.} In fact, this translation emphasizes a key facet of what Poncelet meant by continuity, which is reflected in another expression Poncelet used, when, to refer to the principle, he spoke of the ``principle of the continuity or \textit{permanence of mathematical relations of the figurate magnitude} (principe de la continuit\'e ou \textit{permanence des relations math\'ematiques de la grandeur figur\'ee})'' \cite[p.~319, fourth manuscript, my emphasis]{39}. The intimate connection that the notion of ``continuity'' and that of ``permanence of relationships'' -- as defined by Poncelet -- have in his eyes, is confirmed if we compare this expression to the one mentioned earlier, which reads: ``the \textit{principle} or \textit{the law of continuity} of mathematical relations of the abstract and figurate magnitude (le \textit{principe} ou \textit{la loi de continuit\'e} des relations math\'ematiques de la grandeur abstraite et figur\'ee)''.

In other words, the idea of ``continuity'' here connotes the motion applied to the system that ensures that the validity of the properties and relations \textit{established} for one of its general states \textit{persists} for another state -- this is, the meaning that Cauchy attributes to the term, when he criticizes the principle in his report \cite[p.~72--73]{5}. However, for Poncelet, ``continuity'' simultaneously connotes the persistence of the properties and of the relations between its parts that enter into the \textit{definition} of the figure -- since, as we have seen, the definition of a figure refers to all the positions that it could take ``without \textit{violating the laws}, the \textit{conditions} and the \textit{links} that \textit{subsist} between the various \textit{parts} of the \textit{system}; (\textit{sans violer} les \textit{lois}, les \textit{conditions}, la \textit{liaison} qui \textit{subsistent} entre les diverses \textit{parties} du \textit{syst\`eme})''.\footnote{Del Centina \cite{15} also emphasizes the connection between continuity and permanence (understood in these terms). In an article in preparation, I show that this is precisely the point in which Poncelet's principle of continuity differs from Leibniz's principle with the same name. In this other article, I discuss Poncelet's conception of continuity in analysis and its connection with the continuity in geometry as he understands it. } Whether we speak of the relationships defining the figure or of those obtained by proof, their ``permanence'' does not mean that they stay identical since, as the last part of the latter quotation makes clear, their formulation changes in relation to transformations in the figure.

These assertions can be illustrated using the secants and the chords of a conic section -- the very context in which Poncelet introduced precisely the notion of an ``ideal chord''.\footnote{\cite[p.~26ff]{36}. See \cite{2} for greater detail. Friedelmeyer \cite[p.~118--119]{20} and Lorenat \cite[p.~92--94]{29} comment on the same figure to discuss the introduction of ``ideal chords''.}

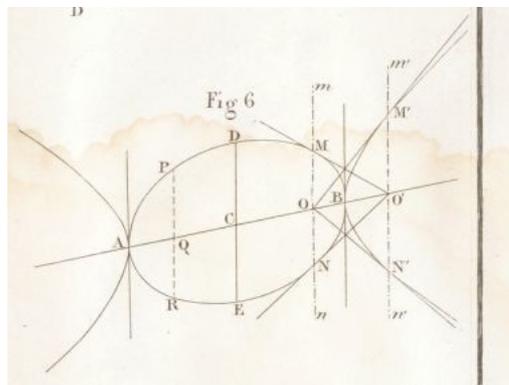
\begin{figure}[t]\centering
 \begin{tikzpicture}[smooth,scale=1,thick,font=\small]
\draw[rotate=26] (0,0) ellipse (1.7 and 1.2);
\draw[-] (-3,-0.6) -- (3.4,0.7);
\draw[-] (-1.62,-1.7) -- (-1.62,1.2);
\draw[-,dashed] (-0.93,-1.25) -- (-0.93,0.86);
\draw[-] (0,-1.25) -- (0,1.25);
\draw[-,dashed] (1.18,-1.5) -- (1.18,2.25);
\draw[-] (1.615,-1.25) -- (1.615,1.8);
\draw[-] (0.4,1.55) -- (2.32,0.48) -- (0.4,-1.38);
\draw[-] (3.5,3.5) -- (1.18,0.25) -- (3.25,-1.9);
\draw[-,dashed] (2.32,-1.5) -- (2.32,2.6);
\draw[-] (3.5, 3.3) .. controls(1.05,0.3) .. (3.25, -1.83);
\draw[-] (-3.25, 1.5) .. controls(-1.145,-0.27) .. (-2.8, -2.2);
\node[scale=0.7] at (-1.8,-0.2) {$A$};
\node[scale=0.7] at (-1.08,0.9) {$P$};
\node[scale=0.7] at (-0.93,-1.4) {$R$};
\node[scale=0.7] at (-0.76,-0.3) {$Q$};
\node[scale=0.7] at (-0.15,0.15) {$C$};
\node[scale=0.7] at (0,1.4) {$D$};
\node[scale=0.7] at (0,-1.4) {$E$};
\node[scale=0.7] at (1.48,0.45) {$B$};
\node[scale=0.7] at (1.03,0.35) {$O$};
\node[scale=0.7] at (1.35,2.1) {$m$};
\node[scale=0.7] at (1.35,1.2) {$M$};
\node[scale=0.7] at (1.35,-0.65) {$N$};
\node[scale=0.7] at (1.35,-1.4) {$n$};
\node[scale=0.7] at (2.5,0.35) {$O'$};
\node[scale=0.7] at (2.5,2.5) {$m'$};
\node[scale=0.7] at (2.5,1.75) {$M'$};
\node[scale=0.7] at (2.5,-0.80) {$N'$};
\node[scale=0.7] at (2.5,-1.4) {$n'$};
 \end{tikzpicture}

\caption{Poncelet's Figure~6, \textit{Trait\'e des propri\'et\'es projectives des Figures}, 1822, \cite[plate~1]{36}.}\label{fig1}
\end{figure}

The line $mn$ being drawn in the plane of a conic section, Figure~\ref{fig1} shows that it intersects the conic section in $M$ and $N$. The midpoint $O$ of the chord $MN$ thereby defined belongs to the diameter $AB$ that is conjugate to the direction of $mn$, as do the centers of all the chords that lines parallel to $mn$ introduce in the conic section (e.g., $DE$, or $PR$). Moreover, the tangents to the conic section at $M$ and $N$ meet in $O'$, which is likewise on the conjugate diameter $AB$. What is more, $A, B, O$ and $O'$ satisfy the following relation:
\[\frac{O'A}{O'B}=\frac{OA}{OB}.
\]

Even though $O$ and $O'$ were constructed with reference to the points $M$ and $N$, the points~$O$ and~$O'$ ``persist'' even if the line $mn$ is continuously moved parallel to itself until it becomes~$m'n'$. The point $O'$ is simply the intersection of $m'n'$ with the diameter $AB$. Moreover, from the point~$O'$, in which the line in its ``new position'' meets with the diameter $AB$, drawing the two tangents to the conic section, $O'M$ and $O'N$ defines a chord $MN$, whose intersection $O$ with $AB$, satisfies the same relation as above with the points $A$, $B$ and $O'$. Poncelet thereby pointed out the ``relation'' between the elements $mn$, $m'n'$, $O$ and $O'$ that holds for the figure and persists wherever the line parallel to $mn$ may be drawn. Given any of them, the other three are determined. Poncelet continued:

\begin{quote}
``since the points $O$ and $O'$, considered as belonging to the straight line $mn$, \textit{subsist} \textit{independently of the reality or the non-reality} of the points of intersection $M$ and $N$ of this straight line and the curve, there is no reason to neglect them in one case rather than in the other; nor is there any reason to neglect the straight line itself, since, when it becomes exterior to the curve, it does not cease to \textit{maintain} certain \textit{dependences} with it. ({\dots}les points $O$ et $O$', consid\'er\'es comme appartenant \`a la droite $mn$, \textit{subsistant ind\'ependamment de la r\'ealit\'e ou de la non-r\'ealit\'e} des intersections $M$, $N$ de cette droite et de la courbe, il n'y a aucune raison de les n\'egliger plut\^{o}t dans un cas que dans l'autre, non plus que cette droite elle-m\^eme, puisqu'en devenant ext\'erieure \`a la courbe elle ne cesse pas de \textit{conserver} certaines \textit{d\'ependances} avec elle.)'' \cite[p.~27, my emphasis]{36}.
\end{quote}

The example illustrates how in any of its positions, a line with the same direction as $mn$ is not independent of the curve but stays in the same relations with it. To put it differently, the same geometrical entities can be seen from different viewpoints. Viewing $O$ as the center of the chord $MN$ might lead us to consider that it vanishes when $M$ and $N$ vanish. However, points such as~$O'$ are tied to the figure by other \textit{relations} (that is, being the intersection of the line~$m'n'$ with the diameter $AB$), and these other relations show that these points~$O'$ \textit{subsist}. What is more, as a result, the part of the system that the line $m'n'$ constitutes has \textit{enduring relationships} with the other part of the system constituted by the curve, even when they do not meet. Since the relationships persist, Poncelet chose to use terms that reflect this ``continuity'' to refer to the geometrical objects. This led him to discard the usual name of $m'n'$ as the ``polar line of $O$'', and to refer to it rather as an ``ideal secant'' of the conic section \cite[p.~27]{36}. Similarly, he referred to $O'$ as the ``ideal center of the ``imaginary chord'' determined by $m'n'$ and to $M'N'$ as an ``ideal chord'', where, for him, $M'$ and $N'$ represent the imaginary points of intersection.\footnote{For a more detailed analysis of the differences between these two types of ideal elements, see \cite{2}. To clarify Figure~\ref{fig1}, when the line $mn$ intersects the conic section in $M$ and $N$, these points satisfy the relation $OM^2=ON^2=p.OA.OB$. The imaginary points of intersection between the conic section and the line in the position $m'n'$ continue to satisfy the relationship $O'M^2~=O'N^2=-p.O'A.O'B$, where the difference of situation implies taking $O'B$ with a negative sign of position (opposed to the sign of $OB$, following Lazare Carnot's way of dealing with signs of position). However, for each point $O'$, for which the line $m'n'$ does not meet the original conic in real points, one may introduce points $M'$ and $N'$ that are such that $O'M'^2=O'N'^2~=p.O'A.O'B$. This relationship defines a conic, to which Poncelet refers as supplementary to the original conic, and which he represents on Figure~\ref{fig1}. Lorenat \cite[p.~94]{29} explains this construction. Poncelet suggests considering $M'N'$ as ``representing, in a fictive manner, the imaginary chord'' corresponding to the secant $m'n'$ \cite[p.~29]{36}. Poncelet refers to the line $M'N'$ as an ``ideal chord''.} We note the correlation between the ``permanence'' of the relationship and the ``permanence'' of the name (e.g., secant, chord or center), Poncelet using the epithet ``ideal'' to qualify in which respect, in some circumstances, the name applies.

Why, one might ask, does Poncelet here distinguish between ``ideal'' and ``imaginary''? We have reached the crux of our argument. Let us read the clearcut explanation that Poncelet develops right away:

\begin{quote}
In general, one uses the adjective \textit{imaginary} to designate any object that, from absolute and real that it was in a given figure, would have \underline{become entirely impossible or inconstruc-} \underline{tible} in the \textit{correlative} figure: ``the one that is supposed to derive from the first figure\footnote{Note of the translator: that is, the primitive figure. This terminology (correlative, primitive\dots ) comes from Lazare Carnot's geometry of position (see~\cite{9}).} by the progressive and continuous motion of some of its parts, without violating the primitive laws~of the system''; the epithet \textit{ideal} would serve to designate the specific mode of existence of an object which, while \underline{remaining}, to the contrary, (KC: that is, in contrast with the previous case) \underline{real in the transformation} of the primitive figure, would nevertheless \underline{cease to depend in an} \underline{absolute and real manner on other objects that define it graphical-} \underline{ly}, because these objects would have become imaginary.\footnote{This last sentence is translated in \cite[p.~94]{29}. However, Lorenat translates ``en demeurant au contraire r\'eel'' as ``in becoming, to the contrary, real'', and not, as I do here, as ``while remaining, to the contrary, real''. } For, just as one already has names in geometry to express the \underline{various modes of existence} that one wants to compare, such as \textit{infinitely small} and \textit{infinitely large}, one must also have names to \underline{express those of} \underline{non-existence}, in order to give accuracy and precision to the language of geometric reasoning.

By comparison with all the definitions that one might want to substitute for them, these definitions have the advantage that they can be \underline{extended} directly to any point, line and surface [{\dots}]; they are useful to shorten the discourse and \underline{extend} the object of geometric conceptions; lastly, they allow us to establish a \underline{contact}, if not always real, \underline{at least fictive, between figures} that at first sight seem to have no relationship with each other, and to \underline{discover easily the relations} \underline{and properties that are common to them}.\footnote{This last paragraph is also translated in \cite[p.~96--97]{29}.}

(En g\'en\'eral on pourrait d\'esigner par l'adjectif \textit{imaginaire} tout objet qui, d'absolu et r\'eel qu'il \'etait dans une certaine figure, serait \underline{devenu enti\`erement impossible ou inconstruc-} \underline{tible} dans la figure \textit{corr\'elative}: `Celle qui est cens\'ee provenir de la premi\`ere par le mouvement progressif et \textit{continu} de quelques parties, sans violer les lois primitives du syst\`eme'; l'\'epith\`ete \textit{id\'eal} servirait \`a d\'esigner le mode particulier d'existence d'un objet qui, en \underline{demeurant} au contraire \underline{r\'eel dans la transformation} de la figure primitive, \underline{cesserait} ce\-pendant de \underline{d\'ependre d'une mani\`ere} \underline{absolue et r\'eelle d'autres objets qui le d\'efinissent} \underline{graphiquement}, parce que ces objets seraient devenus imaginaires. Car, de m\^eme qu'on a~d\'ej\`a en G\'eom\'etrie des noms pour exprimer les \underline{divers modes d'existence} qu'on veut comparer, tels que \textit{infiniment petits}, \textit{infiniment grands}, il faut aussi en avoir pour \underline{exprimer} \underline{ceux de la non-existence}, afin de donner de la justesse et de la pr\'ecision \`a la langue du raisonnement g\'eom\'etrique.

Ces d\'efinitions ont, sur toutes celles qu'on pourrait leur substituer, l'avantage de pouvoir \underline{s'\'etendre} directement \`a des points, des lignes et des surfaces quelconques; ({\dots}) elles servent \`a abr\'eger le discours et \`a \underline{\'etendre} l'objet des conceptions g\'eom\'etriques; enfin elles permettent d'\'etablir un \emph{point de contact}, sinon toujours r\'eel, \emph{au moins fictif, entre des figures} qui paraissent, au premier aspect, n'avoir aucun rapport entre elles, et \emph{de d\'ecouvrir sans peine les relations et les propri\'et\'es qui leur sont communes}.) \cite[p.~28, Poncelet's emphasis is in italics; I underline what I emphasize]{36}.
\end{quote}

In other words, Poncelet distinguished between ``imaginary'' and ``ideal'' as ``two modes'' of ``non-existence''. The former can be illustrated by the points of intersection between the secant $m'n'$ and the conic section, which are imaginary, while the latter is illustrated precisely by the ``ideal secant'' $m'n'$ and the ``ideal chord'' $M'N'$: these entities can still be constructed (since one knows their direction, and one of their points, and also for the chord, its length and the fact that a given point is its middle point), although their dependence on the curve has become imaginary. As the final part of the quotation makes clear, the promotion of these terms derives from Poncelet's reflections on how to use language in science as well as on how to instill generality in geometry. Indeed, the generality arising from the use of these terms derives from the fact that, instead of using different terms for $mn$ and $m'n'$, referring to them as ``secant'' and ``ideal secant'', respectively, allows the geometer to underline their common nature of being ``secant''.\footnote{Lorenat \cite{29} focuses on the issue of language in geometry at the time, comparing ``Louis Gaultier's radical axes, Jean-Victor Poncelet's ideal common chords, and Jakob Steiner's lines of equal power'', and emphasizes the point I am making here (p.~96). She highlights the relation, in Poncelet's eyes, between the new way of naming and the intention of achieving in geometry a generality similar to that of algebraic analysis (p. 98) In particular, Lorenat \cite[p.~93]{29} translates the passage by Poncelet that I summarize here. Note that ``radical axes'' are lines whereas ``common chords'', in contrast with ``secants'', are segments with specific lengths. The objects are related, but distinct. Lorenat's translation of Poncelet's statement in \cite[p.~93]{29} does not distinguish clearly between the two. As a result, she translates Poncelet's statement about the line $m'n'$ ``that these points of intersection with the curve are \textit{imaginary}, and consequently the corresponding chord is itself an \textit{ideal secant} of this curve''. However, it should read: ``that its points of intersection with the curve, and hence the corresponding chord, are \textit{imaginary}, and that it (the line $m'n'$) is itself an \textit{ideal secant} of this curve''. See also \cite[p.~56]{2}.} Note that here, Poncelet employed ``ideal'' as an epithet, to signal when the use of the name in question requires considering imaginary dependences between these elements and real elements of the figure. As a consequence, figures such as Figure~\ref{fig1} need to be read in quite a specific way, and certainly not as a figure in ordinary geometry.\footnote{On further artificial facets of Poncelet's figures, see \cite{2}. Lorenat \cite[p.~176, 182]{28} introduces the difference between Poncelet's concepts of ideal and imaginary and evokes how Poncelet distinguished these two modes of non-existence and modes of existence. She discusses the meaning and part played by geometric evidence for Poncelet. As far as I can tell, her discussion of Poncelet's figures does not dwell on the artificial features of these figures. } Here, the point is to use general enough concepts and general practices with figures and properties.

Moreover, clearly, Poncelet here defined ideal and imaginary elements in full generality, and he insisted that the definitions apply broadly. In the treatise, he was to qualify as ``ideal'' various types of geometrical elements. For instance, seen as the point of intersection of the imaginary tangents at points $M'$ and $N'$, $O$ is an ideal point \cite[p.~30]{36}. However, ``ideal secants'' and ``ideal chords'' remain the key ideal elements that he used. In fact, in the 1820 memoir, even though the notion of ideality is clearly the same, ``ideal secants'' and ``ideal chords'' are the only ``ideal elements'' put into play.\footnote{See the discussion in \cite[p.~367--370]{39}, which is quite similar to what I have just quoted.}

\subsection{Ideal elements as an outcome of Poncelet's philosophical work}\label{section3.3}

The previous quotation is interesting for an additional reason. Indeed, it highlights that Poncelet conceived of the distinction between ideal and imaginary as that between different ``modes of'' ``non-existence'', and that he juxtaposed them in this respect to the ``infinitely small'' and the ``infinitely large'', understood as ``names'' expressing ``various modes of existence''. This way of formulating the opposition is striking because the various manuscripts published by Poncelet show that Poncelet varied significantly in the appreciation of the understanding of and the relationship between these notions. I will limit myself to some remarks on this point, in order to show that Poncelet's 1822 terminology was the outcome of a philosophical reflection, which seems to have reached a stable state in 1822.

The fourth manuscript, published in \cite{39} and completed during the winter 1818--1819, which was devoted to the principle of continuity, contains very important pieces of evidence in this respect. To begin with, clearly for Poncelet, the idea of making use of ``beings of non-existence (\^etres de non-existence)'' in mathematics is due to the practice of analysis. Indeed, he writes:

\begin{quote}
the use and routine of algebraic analysis over the last two centuries\footnote{A few dozens of pages before, in the same notebook, Poncelet \cite[p.~299--300]{39} distinguished between Vieta's application of algebra to geometry and Descartes's, in the way the latter was developed notably by Euler, Lagrange and Monge. I will return to this opposition in the article in preparation on the principle of continuity.} ({\dots}) [have] \textit{introduce}[ed] and consecrat[ed], in the ordinary language of Geometry, certain \textit{expressions} which recall and characterize \textit{beings of non-existence}, and even make us consider these \textit{fictive} beings \textit{as real objects} in reasoning and conception. (L'usage et l'habitude de l'Analyse alg\'ebrique pendant deux si\`ecles ({\dots}) [ont] \textit{introdui}[t] et consacr[\'e], dans le langage ordinaire de la G\'eom\'etrie, certaines \textit{expressions} qui rappellent et caract\'erisent les \textit{\^etres de non-existence}, et font m\^eme consid\'erer ces \^etres \textit{fictifs} comme des \textit{objets r\'eels} dans le raisonnement et la conception.) \cite[p.~344, my emphasis]{39}.
\end{quote}

The example that Poncelet gives to illustrate this statement is highly meaningful for us. It reads as follows:
\begin{quote}
This is how \textit{infinitely large} and \textit{infinitely small} [entities], whose existence is purely hypothetical, have already been introduced into geometry; For, when, e.g., we admit that two parallel straight lines meet at infinity, we \textit{ideally} give an \textit{indefinite} \textit{existence} to their common point of intersection, and thereby also establish the idea of \textit{continuity} in the possible \textit{motion} of this point. On the other hand, the word \textit{infinite} recalls the \textit{geometrical non-existence} of the object designated, and it indicates that the \textit{assumed existence} is purely \textit{ideal}. (C'est ainsi que l'on a d\'ej\`a re\c{c}u en G\'eom\'etrie les \textit{infiniment grands} et les \textit{infiniment petits}, dont l'existence est purement hypoth\'etique; car, lorsqu'on admet, par exemple, que deux droites parall\`eles se rencontrent \`a l'infini, on donne \textit{id\'ealement} une existence ind\'efinie \`a leur point commun d'intersection, et par l\`a aussi on \'etablit l'id\'ee de la \textit{continuit\'e} dans le \emph{mouvement} possible de ce point. D'un autre c\^{o}t\'e, le mot \textit{infini} rappelle la \textit{non-existence g\'eom\'etrique} de l'objet d\'esign\'e, et indique que l'existence suppos\'ee est purement \textit{id\'eale.}) \cite[p.~344--345]{39} (Poncelet emphasizes infinitely large, infinitely small, and infinite. The other italicized words are my emphasis).
\end{quote}

From this quotation, we will retain only three points.
Firstly, the way in which Poncelet used the term ``continuity'' in relation to the point at infinity does seem to combine the idea of the continuity of the ``motion'' as well as that of the continuity of the relationship between lines (namely, that they intersect in one and a single point), thereby confirming what we have seen above.\looseness=-1

Secondly, in contrast to Poncelet's 1822 quotation with which we concluded the previous section, here, in Poncelet's view, the point at infinity, which is the point of intersection of parallel lines, has no ``geometric existence''. In other words, between 1818--1819 and 1822, Poncelet recast his idea about existence and non-existence in geometry.\footnote{The issue of ``existence'' becomes very prominent in Poncelet's reflections starting from the fourth notebook, composed during the winter 1818--1819.}

Thirdly, in the same quotation, Poncelet uses the term ``ideal'' in a way completely different from the use that is evidenced from his 1820 memoir onwards. Indeed, in 1818--1819, the adjective ``ideal'' is applied to the type of existence that is given to a point at infinity. A few pages later, Poncelet has exactly the same use of the adjective ``ideal'', when analyzing the sentence ``these two points have become imaginary'': in his view, the expression is introduced to ``maintain for these points an existence of sign, at least \textit{ideal}, in discourse and conception'' \cite[p.~350, my emphasis]{39}.\footnote{In fact, in the previous manuscript (the third notebook), written between 1815 and 1817, Poncelet used ``ideal'' as interchangeable with ``imaginary'' and ``inconstructible'' and in parallel with infinitely small and large entities. From a given perspective that we need not make explicit, he asserts: ``The same thing will still hold true for the whole interval where these distances remain inconstructible, ideal, imaginary, and the one where they vanish or become infinite momentarily, transitorily (La m\^eme chose aura lieu encore pour tout l'intervalle o\`u ces distances resteront inconstructibles, id\'eales, imaginaires, et celui o\`u elles s'\'evanouissent ou deviennent infinies momentan\'ement, transitoirement.)'' \cite[p.~198]{39}. } Still in the same work, Poncelet employs the same term ``ideal'' to refer to the nature of the relations and properties that are claimed by continuity, when they bring into play elements of the two latter types (infinitely large and infinitely small elements as well as imaginary elements).

Exactly the same use of ``ideal'' as well as the connection between this use and the principle of continuity are also manifest in a letter that he wrote on November 23, 1818 -- that is, more or less at the same time as he was working on this fourth notebook. Indeed, he wrote the following to Olry Terquem, one of the friends with whom, as we saw above, he communicated his mathematical ideas: ``The axiom examined so far, when considered from a certain perspective, is in the end nothing more than the \textit{principle} of permanence or \textit{indefinite continuity of the mathematical laws of magnitudes that are variable by imperceptible succession}, a continuity that, for certain states of the same system, often subsists only in a purely abstract and \textit{ideal} fashion (L'axiome jusqu'ici examin\'e n'est au fond, quand on le consid\`ere sous un certain point de vue, que le \textit{principe} de permanence ou \textit{continuit\'e ind\'efinie des lois math\'ematiques des grandeurs variables par succession insensible}, continuit\'e qui pour certains \'etats d'un m\^eme syst\`eme ne subsiste souvent que d'une mani\`ere purement abstraite et \textit{id\'eale}.)'' \cite[p.~533, except for ``ideal'', Poncelet's emphasis]{39}.

To conclude, the term ``ideal'', which is consistently used in this way in the fourth notebook, changes meaning in Poncelet's writings between 1818--1819 and 1820. In 1818--1819, ``ideal'' qualifies the way in which non-existent elements are nevertheless given a certain type of existence, and it also qualifies the way in which relations and properties extend to them. In contrast, starting from 1820, ``ideal'' refers to ``elements'' that are real, but whose connection to the other parts of the figure has become imaginary \cite[p.~368]{39}. Moreover, this change in the meaning of the term ``ideal'' is, as we have seen, correlated with the change in the understanding of the nature of the infinite magnitudes or elements at infinity.

What is striking is that, still in the same -- fourth -- notebook from the winter of 1818--1819, Poncelet formulated a research program, as follows:

\begin{quote}
({\dots}) \underline{just as one has} \underline{\textit{names}} \underline{to express the} \underline{\textit{various}} \underline{\textit{modes of existence}} \underline{that one wants} \underline{to \textit{compare}}, one must \underline{also have} \underline{\textit{names to express the various modes of non-existence}}, \underline{in} \underline{order to give both accuracy and precision to the language of geometrical reasoning}. Fi- \linebreak nally, \textit{metaphysical notions} themselves have their true \textit{source} in the \textit{perseverance} to apply the idea of \textit{continuity} to a purely \textit{ideal} state of a system, and to \textit{extend to it conceptions and laws that belonged only to the primitive and real state of this system}.

From these ideas, we can see how we should go about tracing a faithful and complete \textit{table of all the metaphysical notions} that belong to \textit{figurate magnitude}. My intention is not to undertake this [table]. It would require more boldness and talent than I possess to dare deliver such a work to the judgment of the severe criticism of geometers: it will suffice for me to have indicated to those more skillful than I am the road that must be followed, and to have taken the first steps along it.

(({\dots}) \underline{de m\^eme qu'on a des} \underline{\textit{noms}} \underline{pour exprimer les} \underline{\textit{divers}} \underline{\textit{modes d'existence}} \underline{qu'on veut} \underline{\textit{comparer}}, il faut aussi \underline{en avoir pour} \underline{\textit{exprimer ceux de la non-existence}}, \underline{afin de donner} \underline{\`a~la fois} \underline{de la justesse et de la pr\'ecision \`a la langue du raisonnement g\'eom\'etrique.} Enfin, les \textit{notions m\'eta\-physiques} elles-m\^emes ont leur \textit{source} v\'eritable dans la \textit{pers\'ev\'erance} qu'on a~d'appliquer l'id\'ee de la \textit{continuit\'e} \`a un \'etat purement \textit{id\'eal} d'un syst\`eme, et d'y \textit{\'etendre les conceptions et les lois qui n'appartenaient qu'\`a l'\'etat primitif et r\'eel de ce syst\`eme}.\looseness=1

D'apr\`es ces id\'ees, on voit comment il faudrait s'y prendre pour parvenir \`a tracer le tableau fid\`ele et complet de \textit{toutes les notions m\'etaphysiques} qui appartiennent \`a la \textit{grandeur figur\'ee}; mon intention n'est pas de l'entreprendre, il faudrait plus de hardiesse et de talent que je n'en poss\`ede pour oser livrer un tel travail au jugement de la critique s\'ev\`ere des g\'eom\`etres: il me suffira d'avoir indiqu\'e \`a de plus habiles que moi la route certaine \`a suivre, et d'y avoir fait les premiers pas.) \cite[p.~345--346, italics refer to my emphasis, whereas the underlining highlights elements that are similar with the 1822 statement quoted above about modes of existence and modes of non-existence]{39}.
\end{quote}

Arguably, the formulation of this research program belongs to what, in the November 23, 1818, letter to Terquem mentioned above, Poncelet referred to as ``the part of [his] work that concerns the metaphysics of geometry (la partie de [s]on travail qui concerne la m\'etaphysique de la g\'eom\'etrie)''.\footnote{\cite[p.~531]{39}. A similar expression (``my efforts to clarify the metaphysics of simple geometry (mes efforts pour \'eclairer la m\'etaphysique de la simple G\'eom\'etrie)'' occurs in Poncelet's response from October 14, 1819 to the two letters in which Terquem communicates Brianchon's as well as Servois's and his own reactions to Poncelet's principle of continuity, which was mentioned above \cite[p.~547]{39}. In the same letter, Poncelet mentioned Brianchon's advice to print separately ``the metaphysical part of [his] work (la partie m\'etaphysique de [s]on travail)'' \cite[p.~545]{39}. The importance that Poncelet attached to mathematicians' research in metaphysics as well as the meaning he gives to that term are topics that are interesting, but not essential to this article. They will thus be the object of another publication.} The way in which this research program in metaphysics is formulated in the quotation above from the winter of 1818--1819 has striking echoes with the long quotation from the 1822 treatise with which we concluded the previous section. Both texts even share the same sentence -- which I have underlined -- on the principles and the epistemological values that should be respected to shape the language to be used in geometry.\footnote{The submitted memoir from 1820 -- fifth notebook -- also testifies to the fact that these ideas had already taken shape in 1820: see, for example, \cite[p.~368, 381--382, 387, 403--404]{39}.} In the text just above, the term ``ideal'' still has its meaning from 1818--1819. By contrast, in the 1822 quotation -- like in the 1820 memoir -- the term takes the new meaning that would remain the one Poncelet would systematically use after 1820. We thus see that the shift of meaning of the term ``ideal'' and the new understanding of the modes of non-existence as well as of those of existence derive from Poncelet's pursuit of the inquiry into ``metaphysics'', for which he was formulating the main questions in 1818--1819. We also see the central part played by continuity understood with the two dimensions emphasized above: continuity in the transformation as well as continuity -- that is, the permanence -- of the relationships.\looseness=1

To conclude, when Kummer discarded the term ``imaginary'' and opted instead for the term ``ideal'', he was making a distinction that Poncelet had made before him, as the outcome of the pursuit of his metaphysical inquiry. As we have seen, for Poncelet in 1822, ``imaginary'' referred to elements in a figure that, after a continuous motion of some of its parts, had become impossible to construct, whereas ``ideal'' qualified those elements that could still be constructed, even though some of their relations to other parts of the figure had become imaginary. This naturally leads us to the question of how Kummer appropriated this distinction. However, before considering this issue, we must present the pieces of evidence we have for the reception of Poncelet's geometrical theories in German-speaking countries.\looseness=1

\subsection{The reception of Poncelet's ideal elements in Germany}\label{sction3.4}

As early as the first issue of the \textit{Journal f\"ur die reine und angewandte Mathematik}, in 1826, the editor August Crelle published a short review of Poncelet's \textit{Trait\'e des Propri\'et\'es Projectives des Figures.}\footnote{\cite[p.~96]{13}. The review is anonymous. However, Poncelet \cite[p.~407]{41} attributes it to Crelle. In 1808, Poncelet, who was a first-year student at the Ecole Polytechnique, had the opportunity of making Alexandre von Humboldt's (1769--1859) acquaintance \cite[p.~406]{41}. He attributes to the latter the fact that Crelle's \emph{Journal} opened its pages to the publications of his articles \cite[p.~407]{41}, and, to begin with, the publication of an article of ``applied mathematics'' in 1827, and then that of \cite{37}. Crelle was a personal friend of Poncelet, their acquaintance dating to before 1830 \cite[p.~407]{41}. Friedelmeyer \cite[p.~146--149]{20} addresses facets of the posterity of Poncelet's treatise in German-speaking countries.} The review insists on the fruitfulness of the use of projections in geometry, and on the interest of the theorems that Poncelet obtained. However, it mentions neither the principle of continuity nor ideal elements.

The situation would change radically with the third volume of the journal, published in 1828. First, Poncelet would publish his first geometrical contribution to the journal and present an overview of his theory \cite{37}. There, he outlined his ideas about the principle of continuity and his motivations for introducing it \cite[p.~215]{37}. He further added:

\begin{quote}
``The \underline{consequence of admitting continuity} was the theory of \textit{ideal secants} and \textit{ideal chords}, and of any object which, \underline{without ceasing to exist effectively} in the transformations of a~single figure, has nevertheless ceased to depend in a purely geometrical manner on other objects to which it related, and which defined or constructed it in the primitive figure. (La \underline{cons\'equence de l'admission de la continuit\'e} a \'et\'e la th\'eorie des \textit{s\'ecantes} et des \textit{cordes} \textit{id\'eales} et, de tout objet qui, \underline{sans cesser d'exister effectivement} dans les transformations d'une m\^eme figure, a pourtant cess\'e de d\'ependre d'une mani\`ere purement g\'eom\'etrique, d'autres objets auxquels il se rapportoit, et qui le d\'efinissoient ou le construisoient dans la figure primitive.) \cite[p.~215--216, Poncelet's emphasis is italicized, mine is expressed through underlining terms]{37}.
\end{quote}

The conception of ideal elements presented here is exactly the same as that to which Poncelet had adhered since 1820. Moreover, the conception is presented in a fully general manner, and not only for ideal secants and chords. One might nevertheless notice that Poncelet did not introduce it in contrast with imaginary elements, as he had done in his treatise.

The outline of Poncelet's theory was formulated in French. However, in the same 1828 issue, Carl Gustav Jacobi (1804--1851) published an article in German, in which he referred to Poncelet's book as ``the famous work \textit{Trait\'e des propri\'et\'es projectives des figures}'' \cite[p.~377]{22}, quoting Poncelet's treatise abundantly. Indeed, Jacobi \cite{22} made use of several theorems that Poncelet established in his book. What is more, in a footnote, Jacobi developed~-- still in German~-- a full explanation for the introduction of ``ideal secants''~-- but not of chords. In fact, Poncelet had made Jacobi's acquaintance in~Paris no later than 1829, and he wrote some recollections of their mathematical conversation on the topic \cite[p.~485--488]{39}. In the same year, Pl\"ucker publishes a monograph in which he referred to ``ideal chords'', by which he meant in fact ``ideal secant'', and which he chose to rename ``Chordal''.\footnote{\cite[p.~49]{34}, quoted in \cite[p.~103]{29}. Lorenat \cite[p.~105]{29} further outlines the reception of Poncelet's 1822 treatise in German-speaking circles. On Pl\textrm{\"u}cker's use of common chords/secants, see \cite[p.~167--175]{27}.}

In the years following these publications, several German works on geometry mention ``ideal secants'' and ``ideal chords''. For instance, in 1830, Christoph Gudermann (1798--1852) published his \textit{Grundriss der analytischen Sphaerik}. In it, he mentioned an opposition between real and ideal chords (for him, \textit{Durchschnittsehne}) as if it were common knowledge, without any specific explanation \cite[p.~138]{21}. Another example is a book related to the lectures in descriptive geometry that Guido Schreiber (1799--1871) gave in the Polytechnic School of Karlsruhe. In it, the author introduced ideal chords (for him, \textit{Durchschnittslinie}), with an explanation and a reference to Poncelet \cite[p.~200, {\S}570]{42}. Again, later, in 1844, Franz Seydewitz (1807--1852), an \textit{Oberlehrer} in the \textit{Gymnasium} of Heiligenstadt, published an article in which he made extensive use of ideal elements -- also in the context of double-column presentation of dual propositions -- making in conclusion an indirect reference to Poncelet \cite[p.~258]{43}.

It is thus not surprising that Kummer knew about ``ideal chords'' -- perhaps through reading German works making use of them, or perhaps through discussions with Jacobi. Why, however, did he choose to illustrate the notion of ideality in geometry using ``ideal chords''? We return to this point in the next section.

However, more importantly for us to understand what meaning Kummer attached to ``ideal chords'', and what inspiration he drew from them, we now need to turn to Michel Chasles's reconceptualization of Poncelet's ideas about continuity and ideality.

\section[Michel Chasles's reflections on contingent and permanent properties]{Michel Chasles's reflections on contingent\\ and permanent properties}\label{section4}

Indeed, in 1837, Michel Chasles published a historical overview of the development of methods in geometry. In this essay, he reviewed the history of geometry from Euclid's \textit{Elements} to the new geometry of his time, his main goal being to highlight how the methods newly acquired by geometry could now compete with the powerful methods that had hitherto been the preserve of analysis applied to geometry, especially from the viewpoint of generality \cite{6, 9}.

In this context, Chasles addressed precisely the type of proof that made use of what Poncelet had introduced under the name of the ``principle of continuity'', referring to it as ``a new mode of demonstration'' \cite[p.~197]{6}. Like Poncelet \cite[p.~309--311]{39}, and with reference to exactly the same example, Chasles attributed to Monge the use of such proofs in mathematical practice \cite[p.~197--198]{6}. Chasles further made it clear that, in contrast to other disciples of Monge's, who followed his practice without justifying it, Poncelet was the first one who made the underlying principle explicit, even though he failed to prove it rigorously \cite[p.~199--200]{6}. Chasles evoked and appeared to embrace Cauchy's criticism thereof (which we have noted in Sections~\ref{section3.1} and~\ref{section3.2}).

Accordingly, Chasles offered an alternative analysis of the principle underlying Monge's practice. The formulation that he gave right away for this principle made no reference to ``continuity'' \cite[p.~198]{6}. A few pages later, he clarified that he preferred not to make use of the notion of ``continuity'', arguing that this would bring in the ``idea of the infinite''.\footnote{On how Chasles deals with elements at infinity, see \cite[p.~746]{6}. } Nor did Chasles refer to any ``ideal'' elements. Instead, he formulated his principle, using new concepts, some of which feature in the name he chooses for this principle: ``the principle of contingent relationships''~\cite{6}. Since this reformulation and Chasles's related reflections will play a key role in Kummer's introduction of ideal elements, I will now examine them in detail.

Like Poncelet, Chasles considers figures as made of parts that could be, e.g., ``points, planes, lines or surfaces'' \cite[p.~198]{6}. On this basis, he begins by introducing the key notion of ``circumstances of construction'' \cite[p.~199--200]{6}. This notion allows Chasles to distinguish, for the same figure, between various dispositions that these elements can have with respect to one another. Some of these dispositions derive from ``particular circumstances of construction'' \cite[p.~199]{6}. If we consider the example of a figure consisting of two circles in a plane, this is the case when these circles are tangent to each other. However, other dispositions derive from ``general circumstances of construction'' \cite[p.~200]{6}. The latter can be illustrated by the case when the two circles cut each other, or by that in which they do not cut each other. In other words, a figure can present different cases, which, Chasles emphasizes, have the ``same generality'' \cite[p.~199]{6}. To analyze the relationship between the same figure drawn in different general circumstances of construction, Chasles introduces the following fundamental opposition:

\begin{quote}
 a figure can present \underline{in its most general construction} two \underline{cases}; in the first one, some \underline{parts} (points, planes, lines, or surfaces), \underline{on which the general construction} of the figure \underline{does not} \underline{necessarily depend,} but that are \textit{\underline{contingent}} or \underline{accidental consequences} of it, are \underline{real and palpable}; in the second case, these \underline{same parts no longer appear}; they have \underline{become imaginary}; and yet the \underline{general conditions} of \underline{construction} of the figure remained the same. (Une figure peut pr\'esenter dans \underline{sa construction la plus g\'en\'erale} deux \underline{cas} diff\'erens (sic); dans le premier, certaines \underline{parties} (points, plans, lignes ou surfaces) \underline{d'o\`u ne d\'epend pas n\'ecessairement la construction g\'en\'erale} de la figure, mais qui en sont des cons\'equences \textit{contingentes} ou accidentelles, sont \underline{r\'eelles et palpables}; dans le second cas, ces \underline{m\^emes parties n'apparaissent plus}; elles sont \underline{devenues imaginaires}; et cependant les \underline{conditions g\'en\'erales} de \underline{construction} de la figure, sont rest\'ees les \underline{m\^emes}.) \cite[p.~198, italicized terms are Chasles's emphasis, whereas I expressed my emphasis using under\-lining]{6}.
 \end{quote}

Some parts in the figure are ``real'' in one case, whereas they ``become imaginary'' in the other. These are, for instance, the two points of intersection in the plane figure consisting of two circles. Chasles refers to such parts as ``contingent or accidental'', in opposition to the parts of the figure which he refers to as the ``\textit{integral} \textit{and} \textit{permanent} \textit{parts} of the figure'' \cite[p.~200]{6}. The latter parts ``belong to its general construction'', independently of the case considered, and ``they are always real'' \cite[p.~200]{6}.

Even though the distinction between the various ``general circumstances of construction'' is important, Chasles insists, it vanishes in the context of ``the application of finite analysis to geometry (l'application de l'analyse finie \`a la g\'eom\'etrie)'', (\textit{ibid.}), which does not discriminate between these cases when establishing ``theorems concerning the \textit{integral} \textit{and} \textit{permanent} \textit{parts} of the figure'' (\textit{ibid.}, Chasles's emphasis). Such results consequently apply to all cases, regardless of whether the ``contingent parts'' are real or imaginary \cite[p.~200, his emphasis]{6}. For Chasles, reasonings of this kind, which ``rely on the general procedures of analysis'' justify -- with some caveats -- \textit{a posteriori} the ``principle of contingent relationships'', which underlies Monge's method and which Chasles formulated as follows:

\begin{quote}
To define this method, we will say that it consists in considering the figure about which we have to prove some general property, in general circumstances of construction, where the presence of certain points, planes or lines, which in other circumstances would be imaginary, makes the demonstration easier. We then apply the theorem that we have proved in this way to the cases of the figure where these points, planes and lines would be imaginary; that is, we consider the theorem as true in all the general circumstances of construction that the figure to which it refers may present. (Pour d\'efinir cette m\'ethode, nous dirons qu'elle consiste \`a~consid\'erer la figure, sur laquelle on a \`a~d\'emontrer quelque propri\'et\'e g\'en\'erale, dans des circonstances de construction g\'en\'erale [Note: elsewhere, Chasles spoke of ``circonstances g\'en\'erales de construction''], o\`u la pr\'esence de certains points, de certains plans ou de certaines lignes, qui dans d'autres circonstances seraient imaginaires, facilite la d\'emonstration. Ensuite, on applique le th\'eor\`eme qu'on a ainsi d\'emontr\'e aux cas de la figure o\`u ces points, ces plans et ces droites seraient imaginaires; c'est-\`a-dire, qu'on le regarde comme vrai dans toutes les circonstances de constructions g\'en\'erales que peut pr\'esenter la figure \`a laquelle il se rapporte) \cite[p.~198]{6}.
\end{quote}

In other words, in Chasles's terminology, Monge's method consisted in proving a theorem about a figure by relying on one of its states drawn in general circumstances of construction. The theorems considered bear only on permanent parts of the figure. However, for such theorems, the practitioner would choose, between these different cases, one that would make the proof easier, owing to the fact that some contingent parts would be real. The truth of the theorem would actually be established for this case. However, once this proof is carried out, the theorem could be claimed to be valid indiscriminately for the figure drawn in any of its general circumstances of construction.\footnote{Chasles \cite[p.~203 and 205]{6} gives alternative formulations of the principle.}

Chasles did not stop here. He considered this practice as unsatisfactory insofar as its validity rests on using analytical reasoning \textit{a posteriori}. In order to remain closer to the ``rigorous principles of the Ancients'', Chasles \cite[p.~205]{6} therefore invited the practitioner to go further and, he said, deeper. Indeed, once a result is established in this way, he suggested that one might look for another proof as follows:

\begin{quote}
Indeed, after having discovered a truth using Monge's somewhat superficial method, which \textit{grasps and makes use} of some \textit{external and palpable circumstance}, which, however, is \textit{fortuitous and ephemeral}, one will have (in order to \textit{establish this truth on reasons that are permanent and independent of variable circumstances of construction} of the figure) to get to the bottom of things and not make use any longer, like Monge, of \textit{secondary and contingent properties}, which, \textit{in some cases, suffice to define} various parts of the figure, but \textit{rather of properties of these same parts of the figure} that are \textit{intrinsic and permanent}. (Il faudra en effet, apr\`es avoir d\'ecouvert quelque v\'erit\'e par la m\'ethode, en quelque sorte superficielle, de Monge, qui \textit{s'empare et tire parti de quelque circonstance externe et palpable}, mais \textit{accidentelle et fugitive}, il faudra, dis-je, pour \textit{\'etablir cette v\'erit\'e sur des raisons permanentes et ind\'ependantes des circonstances variables de construction} de la figure, aller au fond des choses et faire usage \textit{non plus} comme Monge, des \textit{propri\'et\'es secondaires et contingentes} qui suffisent, dans certains cas, pour d\'efinir diverses parties de la figure, mais bien des \textit{propri\'et\'es intrins\`eques et permanentes de ces m\^emes parties de la figure}.) \cite[p.~205, emphasis is mine]{6}.
\end{quote}

 The description of an alternative proof that would rely only on ``permanent reasons'' requires a reflection on what counts as ``permanent parts of a figure''. In this context, Chasles introduced a new element that will be crucial for my argument. Chasles pointed out the fact that we might define the \textit{same parts} of a figure in two ways, essentially distinct from the perspective of the proof sought for. The same parts might be ``defined'' using ``contingent properties''. But they might also enjoy ``intrinsic and permanent properties''. Chasles added specification on these properties:\looseness=-1

\begin{quote}
By \underline{intrinsic and permanent properties}, we mean those that \underline{in all cases} could be used to~\underline{define} and \underline{construct parts of the figure} that we have called \textit{integral} or \textit{principal}, whereas the \textit{secondary} and \textit{contingent} properties are those that can disappear and \underline{become ima-} \underline{ginary} in some circumstances of construction of the figure. (Nous entendons par \underline{propri\'et\'es} \underline{intrins\`eques} et \underline{permanentes} celles qui serviraient, \underline{dans tous les cas}, \`a la \underline{d\'efinition} et \`a la \underline{construction des} \underline{parties de la figure} que nous avons appel\'ees \textit{int\'egrantes} ou \textit{principales}; tandis que les propri\'et\'es \textit{secondaires} et \textit{contingentes} sont celles qui peuvent dispara\^{i}tre et \underline{devenir imaginaires} dans certaines circonstances de construction de la figure.) (\cite[p.~205]{6}, Chasles's emphasis is in italics, mine is in the underlined passages.)
\end{quote}

In other words, some parts of the figure might appear to be contingent, because we have chosen a contingent property to ``define and construct them''. However, one might as well identify a permanent property that they fulfil, and that might highlight the fact that these parts are actually permanent. This remark implies that viewing an object as permanent is the outcome of mathematical work, and is not immediately given. For the second type of proof that Chasles has in mind, a task of this kind is indeed essential.

It is in the context of this discussion, about the change of definition of an object from a~contingent property to a permanent one, that Chasles chose to illustrate his point using the figure of two circles in a plane -- that is, precisely the example that Kummer was to mention in his 1846 article. Chasles's formulation of this example reads as follows:

\begin{quote}
``The theory of circles drawn on a plane gives us an example of this distinction that we make between \textit{accidental} properties and \textit{permanent} properties of a figure. The system of two circles always includes the existence of a certain line, whose consideration is very useful for the whole theory. \underline{When the two circles intersect}, this line is their \textit{common chord}, and this \underline{very circumstance suffices} to \underline{define} and \underline{construct} it; here is what we call \underline{one of its} \textit{contingent} or \textit{accidental} \underline{properties}. But \underline{when the two circles do not intersect}, this \underline{property disappears}, even though, nevertheless, the \underline{line still exists} and that its consideration is still \underline{most useful in the} \underline{theory of circles}. One must therefore \underline{define} and \underline{construct} this line \underline{using another} of its \underline{properties}, which would \underline{hold in all the cases of general con-} \underline{struction} of the figure, which is the system of two circles. This will be \underline{one of its} \textit{permanent} \underline{properties} [{\dots} here Chasles began a discussion of the various permanent properties that can be used, starting from that of considering the line as the radical axis of the two circles]. (La th\'eorie des cercles trac\'es sur un plan nous offre un exemple de cette distinction que nous faisons entre les propri\'et\'es \textit{accidentelles}, et les propri\'et\'es \textit{permanentes} d'une figure. Le syst\`eme de deux cercles comporte toujours l'existence d'une certaine droite, dont la consid\'eration est fort utile dans toute cette th\'eorie. \underline{Quand les deux cercles se coupent}, cette droite est leur \textit{corde commune}, et cette \underline{seule circonstance} suffit pour la \underline{d\'efinir} et la \underline{construire}; voil\`a ce que nous appelons \underline{une de ses} \underline{propri\'et\'es} \textit{contingentes} ou \textit{accidentelles}. Mais \underline{quand les deux cercles ne se coupent pas}, cette \underline{propri\'et\'e dispara\^{i}t} quoique la \underline{droite pourtant} \underline{existe toujours}, et que sa consid\'eration soit \underline{encore extr\^emement utile} \underline{dans la th\'eorie des cercles}. Il faut donc \underline{d\'efinir} cette droite et la \underline{construire par quelqu'une} \underline{de ses autres propri\'et\'es}, qui \underline{ait lieu dans tous les cas de construction g\'en\'erale} de la figure, qui est le syst\`eme des deux cercles. Ce sera \underline{une de ses propri\'et\'es} \textit{permanentes} ({\dots})'' (\cite[p.~205--206]{6}, Chasles's emphasis is in italics, mine is in the underlined passages).
\end{quote}

Interestingly, in his \textit{Aper\c{c}u Historique}, Chasles referred to this line as a ``common chord'', although for Poncelet, it would be a ``secant'' (see also \cite[p.~207]{6}). Yet, the very use of the term ``common chord'' indicates that he had Poncelet's terminology in mind. Setting aside the fact that Chasles referred neither to ``continuity'' nor ``ideal'', with the example of the ``common secant'', we can establish a dictionary between Poncelet's terms and Chasles's. This will be useful when we turn to Kummer's ``ideal divisors''.

Indeed, the parts of a system that for Poncelet were imaginary elements were, for Chasles, contingent parts of the figure.\footnote{For Chasles's interpretation of imaginaries, see \cite[p.~207]{6}. } Moreover, Poncelet's real parts would be either permanent or contingent. However, Chasles had a different approach to the geometrical objects which Poncelet would refer to as ``ideal elements''. For him, they were permanent elements for which only a~contingent definition had been put forward. Accordingly, a permanent property should be found for them. This is precisely what we will see Kummer doing in order to define the prime divisors of complex numbers in a fully general and uniform way.

\section{How Kummer's ideal divisors resemble ideal chords}\label{section5}

\subsection{Kummer reader of Chasles' Aper\c{c}u historique}\label{section5.1}

How, one may ask, did Kummer have access to Chasles's approach to these issues? In fact, as early as 1839 -- that is, only two years after the publication of the book -- Chasles's \textit{Aper\c{c}u historique} was translated into German by Ludwig Adolph Sohncke (1807--1853), one of Jacobi's students.\footnote{\cite[p.~62]{31}.} The German translation of the page about circles just quoted reads as follows:\looseness=1

\begin{quote}
Die Theorie der Kreise in einer Ebene bietet uns ein Beispiel von diesem Unterschied dar, welchen wir zwischen den \textit{zuf\"alligen} und den \textit{bleibenden} \underline{Eigenschaften} einer Figur gemacht haben. Das System zweier Kreise gestattet immer das \underline{Vorhandensein} einer gewissen geraden Linie, deren Betrachtung in \underline{dieser ganzen Theorie von Nutzen ist}. Wenn beide Kreise sich schneiden, so ist diese Gerade ihre \textit{gemeinschaftliche Sehne}, und dieser einzige Umstand reicht hin, um sie zu \underline{erkl\"aren}\footnote{Here, Sohncke translated ``d\'efinir'' using ``erkl\"aren (explain,)'' rather than ``definiren (define)'', as he did elsewhere. In the context in which he translates ``d\'efinir/d\'efinition'' using ``erkl\"aren/Erkl\"arung'', the definition refers to something not clearly defined such as ``analysis and synthesis'' \cite[p.~2, 9]{7} or ``porism'' \cite[p.~10]{7}. He also uses this translation when Chasles spoke of ``defining a method'' \cite[p.~198]{6}; \cite[p.~193]{7}). Elsewhere, Sohncke uses ``erkl\"aren'' to translate ``expliquer'' \cite[p.~101]{7} or ``\'eclairer'' \cite[p.~111, 132]{7}, and ``Erkl\"arung'' to translate ``explication'' \cite[p.~148, 149]{7}. Against this background, I think we should expect here ``definiren''. In exactly the same circumstances, Kummer uses ``Definition'' (see the quotation that follows). I thank one of the referees for having pointed out the fact that according to Felix Mueller's \textit{Mathematisches Vokabularium Franz\"osisch-Deutsch und Deutsch-Franz\"osisch} (Leipzig: Teubner, 1900, p.~185), ``erkl\"aren'' could translate ``d\'efinir''. An examination of Sohncke's translation clarifies in which context he used this translation, and this observation fits with Mueller's specification about this translation.} und zu \underline{construiren}; dieses ist eine von den Eigenschaften, welche wir \textit{zuf\"allige} (KC: a single adjective translating \textit{contingent} or \textit{accidental}) genannt haben. Wenn aber die beiden Kreise sich nicht schneiden, so verschwindet diese Eigenschaft, obgleich die Gerade dennoch \underline{immer besteht} und ihre Betrachtung \underline{von ausserordentlichem Nutzen in der Theorie} des Kreises ist. Man muss daher diese Gerade \underline{definiren} und \underline{construiren} \underline{durch irgend eine andre (sic)} \underline{ihrer Eigenschaften}, welche \underline{stattfindet in allen F\"allen} der \underline{allgemeinen Construction} der Figur, welche hier das System der beiden Kreise ist. Diese wird eine ihrer \textit{permanenten} Eigenschaften sein'' \cite[p.~201, italics are Sohncke's, I underline to express my emphasis]{7}.
\end{quote}

Sohncke's translation is excellent, as can be assessed from the small divergences between the German text and the original one -- quoted above -- that I have indicated between brackets. If, now, we compare the German translation of Chasles's discussion of the figure of the two circles with what Kummer said when he makes his analogy between his introduction of ideal divisors and what happens in geometry, a striking phenomenon appears. To highlight the phenomenon, let us repeat the quotation of Kummer's comparison that I have translated above:

\begin{quote}
 To achieve now a lasting \underline{definition} of the true (usually ideal) prime factors of the complex numbers, it was necessary to find out the \underline{properties} of the prime factors of complex numbers that would \underline{persist/remain in all circumstances}, which would be absolutely independent of the \underline{contingency/accidental} circumstances of whether the actual decomposition takes place or not, more or less precisely as, when in geometry, one speaks of the \underline{chord common to two circles} also \underline{when the circles do not intersect each other}, one looks for an actual \underline{definition} of this ideal \underline{common chord} that fits for \underline{all} situations of the \underline{circles}.

(Um nun zu einer festen \underline{Definition} der wahren (gew\"ohnlich idealen) Primfactoren der complexen Zahlen zu gelangen, war es n\"othig, die unter \underline{allen Umst\"anden bleibenden Eigen-} \underline{schaften} der Primfactoren complexer Zahlen aufzusuchen, welche von der \underline{Zuf\"alligkeit}, ob die wirkliche Zerlegung Statt habe, oder nicht, ganz unabh\"angig w\"aren: ohnegef\"ahr ebenso wie wenn in der Geometrie von der \underline{gemeinschaftlichen Sehne zweier Kreise} gesprochen wird, auch wenn die Kreise \underline{sich nicht schneiden}, eine wirkliche \underline{Definition} dieser idealen \underline{gemeinschaftlichen Sehne} gesucht wird, welche f\"ur \underline{alle} Lagen \underline{der Kreise} passt.) \cite[p.~88, underlining is mine]{23}.
\end{quote}

I have underlined all the terms that appear in both contexts. The comparison shows clearly how many terms Kummer has taken from the page that Chasles inserted in his \textit{Aper\c{c}u historique}. Before I explain the importance of this remark, let me begin by establishing why we can be sure that Kummer did read Chasles's \textit{Aper\c{c}u historique}.

 Two hints clearly support this conclusion. The first one lies in the fact that, to refer to the secant line, Kummer, like Chasles, used the term ``common chord'', translated by Sohncke as ``\textit{gemeinschaftlichen Sehne}''. We have seen that other German authors used different terms for the ``common chord'', and hence, for this term too, Kummer uses Sohncke's terminology. However, more importantly, as I have pointed out above, here Chasles used the term ``chord'', although he actually spoke of a common ``secant'' -- the whole line that, as was mentioned, constitutes the radical axis. We note exactly the same phenomenon in Kummer's comparison.\looseness=-1

The second hint derives from a letter that Kummer wrote to Kronecker on 14 June, 1846 \cite[p.~98--99]{46}. This is the second extant letter in which Kummer commented upon the ideal factors and in which, as I have mentioned in the introduction, he put forward the comparison between the composition of numbers using ideal factors and the composition of chemical substances. Incidentally, he mentioned having discussed his ideal divisors with Jacobi. At the end of this letter, however, he turned to a completely different topic, as follows:

\begin{quote}
From my own work, I further mention something really pretty. Indeed, I have solved in a fully general way the problem of finding quadrilaterals, whose sides and diagonals are rational. This problem becomes interesting, insofar as CHASLES (wrongly) thinks that the Indians had set themselves this problem and had solved it in their way and further insofar as the punch line of the whole thing comes down to the fact of making a root $\sqrt{a+bx+cx^2+dx^3+{ex}^4}$ rational, which EULER had already treated in the Algebra and about which JACOBI has later shown that the theory of elliptic function solves this problem (Von meinen eigenen Arbeiten erw\"ahne ich noch eine recht nette Sache. Ich habe n\"amlich ganz allgemein die Aufgabe gel\"ost Vierecke zu finden, deren Seiten und Diagonalen rational sind. Dieselbe wird dadurch interessant da{\ss} CHASLES (f\"alschlich) meint die Inder h\"atten sich diese Aufgabe gestellt und in ihrer Weise gel\"ost, ferner dadurch da{\ss} die Pointe des ganzen darauf hinauskommt eine Wurzel $\sqrt{a+bx+cx^2+dx^3+{ex}^4}~$ rational zu machen, welches EULER schon in der Algebra behandelt und wovon JACOBI sp\"ater gezeigt hat da{\ss} die Theorie der elliptischen Functionen diese Aufgabe l\"ost.) \cite[p.~99]{46}.\looseness=-1
\end{quote}

Kummer's allusion refers to the Note XII of the \textit{Aper\c{c}u historique}, titled ``Sur la G\'eom\'etrie des Indiens, des Arabes, des Latins et des Occidentaux au moyen \^{a}ge'' \cite[p.~416--542]{6}.\footnote{Ivahn Smadja has analyzed Kummer's work on these Sanskrit writings in ``Part 4: Les quadrilat\textrm{\`e}res de Brahmagupta: Chasles, Kummer, Hankel'' \cite[p.~17--31]{44}. See his publication in preparation on this topic.} The letter, dated from a few months after the presentation to the Berlin academy, thus shows that in the same months that he was working on his theory of complex numbers, Kummer was reading Chasles's \textit{Aper\c{c}u historique} in great detail.

\subsection{Kummer's search for a permanent definition of prime factors}\label{section5.2}

Let us return to the terms that, as I have underlined above, in 1846 Kummer \cite{23} borrowed from Chasles's page on the system of two circles. These terms include: ``contingency (Zuf\"alligkeit)'', ``permanent property (bleibende Eigenschaft)'', ``common chord to two circles (gemeinschaftliche Sehne zweier Kreise)'', ``when the two circles do not intersect each other (Wenn aber die beiden Kreise sich nicht schneiden)'', ``definition/to define (Definition/definiren)'', ``Umstand (circumstance)''.

These terms relate to the issue of ``defining'' permanent parts of a system using permanent/contingent properties. We have seen that Kummer put forward his third comparison -- the one that pointed out a parallel between ``ideal divisors'' and ``ideal chords'' -- precisely when he turned to the issue of the definition of divisors. Focusing on these terms highlight a crucial phenomenon: Indeed, these terms are correlated with the very structure of Kummer's first presentation of his ideal divisors in 1846.

On this basis, I will argue that Kummer's 1846 article reflects an important part of the thought-processes that led Kummer to the introduction of his ideal factors. More precisely, Kummer's achievement was inspired by a transfer, to the context of the divisibility of complex numbers by prime factors, of Chasles's prescription about looking for permanent properties to define objects that persist in all circumstances.

 Indeed, immediately after the comparison on which I have dwelled, Kummer wrote:

\begin{quote}
Among the \textit{permanent properties} of the complex numbers which are handy to be \textit{used as definitions} of ideal prime factors, there are \textit{several} which fundamentally always \textit{lead to the same result} and among which I have \textit{chosen} one as the \textit{simplest} and the \textit{most general}. (Dergleichen \textit{bleibende Eigenschaften} der complexen Zahlen, welche geschickt sind,\footnote{Here ``um'' was added in \cite[p.~320]{25}.} \textit{als Definitionen} der idealen Primfactoren \textit{benutzt} zu werden, giebt (sic) es \textit{mehrere}, welche im Grunde immer \textit{auf dasselbe Resultat f\"uhren} und von denen ich eine als die \textit{einfachste} und \textit{allgemeinste} gew\"ahlt habe.) \cite[p.~88--89, my emphasis]{23}.
\end{quote}

Clearly, Kummer had been working on various ``permanent properties'', with the intention to seek an appropriate ``definition'' of ideal prime factors.\footnote{Kummer's first letter about ideal divisors from October 18, 1845 contains a similar hint \cite[p.~96]{46}, without, however, speaking of ``permanent property''.} As he made clear, his choice among the properties he considered was dictated by two epistemological values that were also crucial for Chasles: ``simplicity'' and ``generality'' \cite{9}. This work on definitions is reflected in the first presentation of the ``theory of complex numbers''. Indeed, what is striking against this backdrop is that Kummer's first public presentation \cite{23} has three parts, each devoted to a definition of ideal prime factors, the first definition and then the second one being in turn discarded for their lack of generality and being eventually replaced by the third one. Let us examine them in turn, to make this point clear.

 The first definition considered by Kummer is introduced immediately after the last statement quoted, and it focuses on the decomposition of a prime integer $p$ into complex prime factors, in the case when $p$ is of the form $m\lambda + 1$ (recall that $\lambda$ is prime and $\alpha$ is an imaginary root of the equation $\alpha^\lambda=1$).\footnote{As was mentioned in Section~\ref{section2}, for Kummer, prime complex numbers are such that, when they divide a product of two numbers, they must divide one of these two factors. The condition ``$p$ is of the form $m\lambda + 1$'' corresponds chronologically to the fact that Kummer first explored the complex numbers $f(\alpha)$, whose norm (that is, $f(\alpha) f(\alpha^2)\cdots f(\alpha^{\lambda-1})$) was a prime integer $p$. See, e.g., the letter to Kronecker from April 10, 1844 (in \cite[p.~83--87]{46}). In this letter, Kummer shows that, if $p$ is the norm of $f(\alpha)$, there exists an integer $\xi$ satisfying the property that $f(\alpha )$ divides $\alpha -\xi$, or, in other terms, $\alpha -\xi \equiv 0$ (mod $f(\alpha)$). Consequently, in these cases, $1+ \xi + \xi ^2+\dots + \xi^{\lambda -1} \equiv 0$ (mod~$p$). Either to $\xi \equiv 1$ (mod $p$), and $p = \lambda $, or this implies that $\xi^\lambda \equiv 1$ (mod~$p$), and hence $\lambda$ divides $(p-1)$. We will see that this particular case is also the simplest one. If we except the complex numbers whose norm is $\lambda$, Kummer showed that prime numbers $p$ that are norms of numbers $f(\alpha )$ are of the form $ m \lambda + 1$. This property derives from the fact that in such cases, there exists an integer $\xi$ satisfying the property that $f(\alpha)$ divides $\alpha -\xi$, or, in other terms, $\alpha -\xi \equiv 0$ (mod $f(\alpha )$).\label{footnote65}} The introduction of the first case reads as follows:

\begin{quote}
If $p$ is a prime number of the form $m\lambda + 1$, \textit{in many cases} it can be represented as a product of the following $\lambda -1$ complex factors: $p = f(\alpha) f(\alpha^{2}) f(\alpha^{3})\cdots f (\alpha^{\lambda-1})$;
\textit{however}, in \textit{cases} when a \textit{decomposition} into \textit{actual} complex prime factors is \textit{not possible}, then \textit{ideal prime factors} must be introduced, to get the \textit{same decomposition}. (Ist $p$ eine Primzahl von der Form $m \lambda+ 1$, so l\"asst es\footnote{Kummer \cite[p.~320]{25} replaces ``es'' by ``sie''.} sich \textit{in vielen F\"allen} als Product von folgenden $\lambda -1$ complexen Factoren darstellen: $p = f(\alpha) f(\alpha^{2}) f(\alpha^{3})\cdots f(\alpha^{\lambda-1})$,\footnote{Kummer \cite[p.~320]{25} changes the comma into a semi-colon.} wo \textit{aber} eine \textit{Zerlegung} in \textit{wirkliche} complexe Primfactoren \textit{nicht m\"oglich} ist, da eben\footnote{Kummer \cite[p.~320]{25} replaces ``, so'' with ``: dann''} sollen die \textit{idealen Primfactoren} eintreten, um \textit{dieselbe} zu leisten.) \cite[p.~89, my emphasis]{23}.
\end{quote}

The problem is clearly set: the decomposition into ``actual prime factors'' occurs in ``many cases'', but not in all. The point of introducing ``ideal prime factors'' is to get ``the same decomposition'' uniformly for any such $p$. We recognize an assumption similar to Poncelet's credo in the ``permanence of laws''. The path followed by Kummer's presentation for this first set of prime numbers $p$ that he considered is then to find out an adequate property that ``actual prime factors'' all possess and that can be used as a definition of ``prime factors'' for all cases of divisors. Kummer therefore went on with the statement of a property of actual factors as follows:\looseness=-1

\begin{quote}
If $f(\alpha)$ is an \textit{actual} complex number and a \textit{prime factor} of $p$, then this complex number has the \textit{property} that, if, instead of the root of the equation $\alpha^\lambda = 1$, one substitutes in $f(\alpha)$ a~certain root of the congruence $\xi^\lambda \equiv 1$, mod~$p$, then $f(\xi) \equiv0$, mod~$p$. Therefore also, if the prime factor~$f(\alpha)$ is contained in [KC: contained as a factor, that is, divides] a~complex number~$\Phi(\alpha)$, thus $\Phi(\xi) \equiv 0$, mod~$p$; and \textit{conversely}: \textit{if} $\Phi(\xi) \equiv 0$, mod~$p$, \textit{and} $p$ \textit{can be decomposed} into $\lambda - 1$ complex prime factors, then $\Phi(\alpha)$ contains [KC: contains as a factor, that is, is divided by] the prime factor $f(\alpha)$. (Ist $f(\alpha)$ eine \textit{wirkliche} complexe Zahl und\footnote{Kummer \cite[p.~320]{25} inserts here ``ein''.} \textit{Primfactor von $p$}, so hat sie die \textit{Eigenschaft}, dass wenn anstatt\footnote{Instead of ``wenn anstatt'', \cite[p.~320]{25} has here ``, wenn statt''.} der Wurzel der Gleichung $\alpha^\lambda = 1$ eine bestimmte Congruenzwurzel von $\xi^\lambda \equiv 1$, mod~$p$, substituirt wird, $f(\xi) \equiv 0$, mod~$p$, ist. Darum auch\footnote{Instead of ``Darum auch'' \cite[p.~320]{25} has here ``Also auch, \dots''.} wenn in einer complexen Zahl $\Phi(\xi)$ der Primfactor $f(\alpha)$ enthalten ist, so\footnote{Kummer \cite[p.~320]{25} deletes ``so''.} wird $\Phi(\xi)\equiv 0$, mod~$p$; und \textit{umgekehrt}: \textit{wenn} $\Phi(\xi) \equiv 0$, mod~$p$, \textit{und} $p$ in $\lambda- 1$ \textit{complexe Primfactoren zerlegbar ist,} enth\"alt~$\Phi(\alpha)$ den Primfactor~$f(\alpha)$.) \cite[p.~89, my emphasis]{23}.
\end{quote}

In other words, Kummer pointed out a property of an actual prime factor $f(\alpha)$ of $p$, the norm of which is~$p$. To understand this property, note that, if $f(\alpha)$ is such a factor of a prime integer~$p$, there exists an integer $\xi$ that is such that $\alpha \equiv \xi$ (mod $f(\alpha)$) and that hence satisfies the equation $\xi^\lambda \equiv 1$ (mod~$p$) (see footnote~\ref{footnote65}). As a result, if a prime factor $f(\alpha)$ of~$p$ can be exhibited, there exists an integer~$\xi$ such that $f(\alpha) \equiv f(\xi) \equiv 0$ (mod\ $f(\alpha)$), and hence $f(\xi) \equiv 0$ (mod~$p$). Consequently, the test of the divisibility of $\Phi(\alpha)$ by the prime factor $f(\alpha)$ associated with $\xi$ is equivalent to the test of the integer $\Phi(\xi)$ by~$p$.\footnote{Precise theorems and proofs for these statements can be found in \cite[p.~89--95]{18}.} Now, Kummer proceeded in exactly the same way as Chasles invited practitioners to proceed, when the French geometer suggested defining the ordinary secant no longer as the straight line going through the points of intersection between the circles, but rather using the same property as that for the ``ideal secant''~-- that is, using the permanent property that they are both the radical axis of the two circles: For Kummer, any non-trivial root of the equation $\xi^\lambda \equiv 1$, mod~$p$, gives rise to a test among integers, which enables one to test whether the complex divisor corresponding to this root~$\xi$ in the case of numbers~$p$ of the form under consideration divides any given complex number. Accordingly, having observed the aforementioned property of an actual prime factor, Kummer next turned this property into a general definition of any prime factors of~$p$, including ``ideal divisors''. Indeed, Kummer's subsequent statements read as follows:

\begin{quote}
Now, the \textit{property} $\Phi(\xi) \equiv 0$, mod~$p$ is one which is such that it in itself \textit{does not depend at all} on the \textit{decomposability} of the number $p$ into $\lambda-1$ prime factors; \textit{it can therefore be used} as a \textit{definition}, through which it is decided that the complex number $\Phi(\alpha)$ contains the ideal prime factor of $p$ that belongs to $\alpha=\xi$ if $\Phi(\xi)\equiv 0$, mod~$p$. \textit{Each} of the $\lambda-1$ \textit{complex prime factors} of~$p$ is in this way \textit{replaced by a condition of congruence}. This \textit{suffices} to show that the complex prime factors, \textit{whether they be actual or exist only ideally}, give to the complex numbers the \textit{same defining character}. (Die \textit{Eigenschaft} $\Phi(\xi)\equiv 0$, mod~$p$, ist nun eine solche, welche f\"ur sich\footnote{Kummer \cite[p.~320]{25} inserts here ``selbst''.} von der \textit{Zerlegbarkeit} der Zahl $p$ in $\lambda-1$ Primfactoren \textit{gar unabh\"angig ist};\footnote{Instead of ``gar unabh\"angig ist'', Kummer \cite[p.~320]{25} has ``gar nicht abhangt (sic) (does not depend at all)''.} sie \textit{kann} daher\footnote{Instead of ``daher'', Kummer \cite[p.~320]{25} has ``demnach''.} \textit{als Definition benutzt werden}, indem bestimmt wird, dass die complexe Zahl $\Phi(\alpha)$ den idealen Primfactor von $p$ enth\"alt, welcher zu $\alpha=\xi$ geh\"ort, wenn $\Phi(\xi)\equiv 0$, mod~$p$, ist. \textit{Jeder} der $\lambda-1$ \textit{complexen Primfactoren} von $p$ wird so \textit{durch eine Congruenzbedingung ersetzt.} Dies \textit{reicht hin}, um zu zeigen, dass die complexen Primfactoren, sie seien \textit{wirklich, oder nur ideal vorhanden}, den complexen Zahlen denselben bestimmten Character imprimiren.\footnote{Instead of ``imprimiren'', \cite[p.~320]{25} has ``ertheilen''. The verb ``imprimiren'' has the nuance of ``imprinting'' a similar structure to the numbers, whereas ``erteilen/ertheilen'' can be translated as ``to grant the numbers the same character''.}) \cite[p.~89, my emphasis]{23}.
\end{quote}

The property identified does not depend on the actual decomposability of~$p$ into actual prime complex numbers, in exactly the same way as the property of a line to be the radical axis does not depend on whether the circles intersect. Kummer thus replaced the exhibition of actual prime complex divisors of~$p$ by this permanent property, taken as a definition of a prime factor that relies on a congruence test. The congruence $\xi^\lambda \equiv 1$, mod~$p$, has $\lambda - 1$ roots that differ from~1. Each of these roots is associated with a prime divisor of~$p$.

The definition put forward is that of the divisibility by a factor, and not that of the factor itself. Kummer had identified a property of prime factors of~$p$ that enables him to detect them in all cases: whether they are actual (like a number $f(\alpha)$) or ideal, they satisfy exactly the same property. If we draw on Chasles's ``principle of contingent relationships'', which Kummer was transposing here, this suggests that for Kummer, these ideal factors are ``integral and permanent parts'' of the complex number that they divide. Only a belief of this kind can justify looking for permanent properties. However, unlike the case of the radical axis, the ideal divisors cannot be ``exhibited''.

One might be tempted to believe that Kummer had now managed to define ideal factors and that his work is thus over. However, he immediately pointed out that this first ``permanent property'', which he has introduced as such, will not solve his problem, since, as he states:

\begin{quote}
However, we do \textit{not use} the congruence conditions in the way that is given here \textit{as definition of the ideal prime factors}, because these conditions would \textit{not be sufficient} to \textit{represent} several \textit{identical/equal ideal prime factors} that would occur in a complex number, and \textit{because} these \textit{conditions} are \textit{too limited} and would give \textit{only} ideal prime factors of the real prime numbers of the form $m \lambda+ 1$. (In der hier gegebenen Weise aber \textit{gebrauchen wir die Congruenzbedingungen nicht} als \textit{Definitionen} der idealen Primfactoren, weil dieselbe\footnote{Instead of ``dieselbe'', \cite[p.~320]{25} has ``diese''.} \textit{nicht hinreichend} sein w\"urden, mehrere gleiche, in einer complexen Zahl vorkommende ideale Primfactoren \textit{vorzustellen}, und \textit{weil} sie, \textit{zu beschr\"ankt}, \textit{nur} ideale Primfactoren der realen Primzahlen von der Form $m \lambda+ 1$ geben w\"urden.\footnote{Here, \cite[p.~320]{25} has a typo, since it reads ``$m \lambda-1$''. \cite[p.~89]{23} has the correct expression ``$m \lambda+ 1$''}) \cite[p.~89, my emphasis]{23}.
\end{quote}

 In effect, in what follows, this first ``permanent property'' is replaced by a second, and then a~third one, which respectively extend it and successively encompass the cases that Kummer notes here as not dealt with if one uses the first definition. The last definition will be fully general. We might thus assume that this structure of the presentation reflects Kummer's statement quoted above, when he asserted that he had found several permanent properties, which lead to introducing the same entities, and that he had chosen the ``most general''. This remark highlights the part played by Chasles' invitation to concentrate on ``permanent properties''. What is more, I argue that the structure of the presentation -- that is, the introduction of Kummer's first definition of ideal complex divisors and its subsequent extensions into two increasingly more general definitions -- also appears to reflect part of the process of Kummer's discoveries with respect to ideal factors.

\subsection{Kummer's discovery as a reflection of the impact of geometry}\label{section5.3}

To understand this point, let us consider how Kummer shapes the second definition -- which eliminates the previous restriction on~$p$ and deals with all prime integers $p$ ($\lambda$~excluded), the third definition then turning to all complex numbers and solving the problem of the multiple factors that can occur among their prime divisors.\footnote{I rely on \cite[p.~107--133]{18}, which I summarize coarsely here. The relevant sections of \cite{18} give a detailed and rigorous presentation.} For this, we will need to introduce a few notations and some mathematical facts.

For any prime integer $p$ different from $\lambda$, we have
\[p^{\lambda -1}\equiv 1 \quad (\text{mod} \ \lambda).
\]

Let us call $f$ the smallest integer for which we have
\[p^f\equiv 1 \quad (\text{mod} \ \lambda )\]
and $e$ the integer such that
\[ef=\lambda -1.\]

For the second definition to which Kummer turns, he will make use of a tool introduced by Carl Friedrich Gauss (1777--1855) in the \textit{Disquisitiones Arithmeticae}: the ``periods''. The~$e$ periods, which are all sums consisting of $f$ terms, are defined using a primitive root $\gamma$ of $(\mathbb{Z}/\lambda \mathbb{Z})^{*}$~-- in other words, a multiplicative generator of $(\mathbb{Z}/\lambda \mathbb{Z})^{*}$~-- as follows:
\begin{gather*}
\eta _0=\alpha^{\gamma^e}+\alpha^{\gamma^{2e}}+\dots +\alpha^{\gamma^{ef}},\\
\eta _1=\alpha^{\gamma^{e+1}}+\alpha^{\gamma^{2e+1}}+\dots +\alpha^{\gamma^{ef+1}},\\
\cdots\cdots\cdots\cdots\cdots\cdots\cdots\cdots\cdots\cdots\cdots\cdots\cdots\cdots\\
\eta_{e-1}=\alpha ^{\gamma ^{e+(e-1)}}+\alpha ^{\gamma ^{2e+(e-1)}}+\dots +\alpha ^{\gamma ^{ef+(e-1)}}.
\end{gather*}

The period ${\eta} _0$ (Kummer notes it $\eta $) is the one containing $\alpha$, and hence it is independent of the primitive root chosen, whereas the other periods are simply permuted with one another if one chooses another primitive root. Several facts will be essential for Kummer's second definition.

 First, if $h(\alpha)$ is a prime complex divisor of a prime integer $p$ such that
$p^f\equiv 1$ (mod~$\lambda$), then, in a way similar to what was done for the previous case, one can find $e$ integers $u_0, u_1, \dots, u_{e-1}$ such that, for each integer $i$ comprised between 0 and $e-1$, we have
${\eta} _i\equiv u_i$ (mod $h(\alpha)$).

Second, to test the divisibility of a complex number $g(\alpha)$ by such a prime complex divisor $h(\alpha)$, a key fact consists in noticing that for any $g(\alpha)$, one can find $f$ functions $\varphi _i $ of the periods (which Kummer notes $\varphi _i(\eta)$)\footnote{Kummer justified this notation for a function of the $e$ periods by the fact that any period ${\eta}_i$ can be expressed as a linear combination, with rational coefficients, of one of them, say $\eta$, and its powers $\eta^2, \dots, \eta^{e-1}$ \cite[p.~406, 411]{26}. One can also think of $\varphi_i(\eta)$ as a function of $\eta, \dots, \eta_{e-1}$. The integers $u_i$ satisfy congruences that have exactly the same form as the equations satisfied by the periods \cite[p.~410--411]{26}. The same remarks thus apply about the notations, replacing equalities by congruences (see below).} such that
\[
g(\alpha)= \varphi _0(\eta)+\alpha \varphi _1(\eta)+\dots +\alpha ^{f-1}\varphi _{f-1}(\eta).\]

As a result,
\begin{gather*}
\varphi _0(\eta)+\alpha \varphi _1(\eta)+\dots +\alpha ^{f-1}\varphi _{f-1}(\eta)\equiv \varphi _0(u)+\alpha \varphi _1(u)+\dots +\alpha ^{f-1}\varphi _{f-1}(u) \quad\! (\text{mod} \ h(\alpha))
\end{gather*}
in which the periods ${\eta} _i$ have been replaced by the integers $u_i$ (using a similar notation for the~$u$'s as for the $\eta $'s).

One can prove that two distinct numbers of the kind $\varphi_0(u)+\alpha \varphi _1(u)+\dots +\alpha ^{f-1}\varphi _{f-1}(u)$ are congruent to each other modulo $h(\alpha)$ if and only if they are equal, and hence any complex number $g(\alpha)$ -- which, as we have seen, can be written as $\varphi _0(\eta)+\alpha \varphi _1(\eta)+\dots +\alpha ^{f-1}\varphi _{f-1}(\eta)$ -- is divisible by the prime complex divisor $h(\alpha)$ if and only if
\begin{gather*}
\varphi _0(u)\equiv 0 \quad (\text{mod} \ p),\\
\varphi _1(u)\equiv 0\quad (\text{mod} \ p),\\
\cdots\cdots\cdots\cdots\cdots\cdots\cdots \\
\varphi _{f-1}(u)\equiv 0\quad (\text{mod} \ p).
\end{gather*}

Consequently, divisibility by a prime complex factor can be tested using only congruences among integers. The key fact is that the $e$ periods ${\eta} _i$ satisfy an equation of degree $e$ with integer coefficients.\footnote{The product of any two periods can be expressed as a linear combination of periods with integral coefficients, and the sum of all periods is $-1$.} Solving it as a congruence modulo the prime number $p$ yields $e$ integers $u_i$, which (by a circular permutation among the $u_i$) allow Kummer to define the $e$ ideal complex divisors of any prime integer $p$ (except the prime integer $\lambda$ which is dealt with separately).

On the basis of these notations and facts, let us return to how Kummer introduced the second definition he considers. He began with an outline of the case in which the decomposition can be carried out, as follows:

\begin{quote}
``\textit{Every prime factor of a complex number} is \textit{always} at the same time also a prime factor of some real prime number $q$, and the \textit{constitution of the ideal prime factors} depends above all on the exponent to which $q$ belongs, for the modulus $\lambda$. Let this one be $f$, so that $q^{f} \equiv 1$ mod $\lambda$, and $\lambda-1= e f$. A \textit{prime number} \textit{$q$ of this kind} can \textit{never} be decomposed into \textit{more than} $e$ complex prime factors, which, \textit{if the decomposition can actually be carried out}, are represented \textit{as linear functions of the e periods} with each $f$ terms. I denote these periods of the roots of the equation $\alpha^\lambda = 1$ as $\eta,\eta_1,\eta_2,\dots,\eta_{e-1}$; and hence in the order in which each becomes the next one, when $\alpha$ is turned into $\alpha^\gamma$, where $\gamma$ is a primitive root of $\lambda$. (\textit{Jeder Primfactor einer complexen Zahl} ist \textit{immer} zugleich auch Primfactor irgend einer realen Primzahl $q$, und \emph{die Beschaffenheit der idealen Primfactoren} ist besonders von dem Exponenten abh\"angig, zu welchem $q$ geh\"ort, f\"ur den Modul $\lambda$, derselbe\footnote{Here, \cite[p.~321]{25} starts a new sentence.} sei $f$, so dass $q^{f} \equiv 1$ mod $\lambda$, und $\lambda-1= e f$. \textit{Eine solche Primzahl} $q$ l\"asst sich \textit{niemals} in \textit{mehr als} $e$ complexe Primfactoren zerlegen, welche, \textit{wenn diese Zerlegung wirklich ausf\"uhrbar ist}, sich als \emph{line\"are}\footnote{In \cite[p.~321]{25}, the word is written ``lineare''.} \emph{Functionen der $e$ Perioden} von je $f$ Gliedern darstellen. Diese Perioden der Wurzeln der Gleichung $\alpha^\lambda = 1$ bezeichne ich durch $\eta,\eta_1,\eta_2,\dots,\eta_{e-1}$; und zwar in der Ordnung, dass jede in die folgende \"ubergeht, wenn $\alpha$ in $\alpha^\gamma$ verwandelt wird, wo $\gamma$ eine primitive Wurzel von $\lambda$ ist.)'' \cite[p.~89--90, my emphasis]{23}.
\end{quote}

 We can notice that the description of the prime complex divisors starts with the fully general case (prime factors of a complex number) and introduces prime integers as a key case to deal with the prime complex divisors. The discussion of prime integers is now fully general, and does not give any special status to the case dealt with previously. In fact, the first extant letter in which, on October~18, 1845, Kummer announced his breakthrough to Kronecker, enables us to establish that previously, Kummer had treated the cases of distinct types of prime integers separately (again setting the prime integer $\lambda$ aside). Indeed, about a letter prior to the October~18 one, he wrote:

\begin{quote}
The notations that I have will already be known to you from the previous letter and otherwise, that is, that $\lambda - 1 = e f$ and $q$ belongs to the exponent $f$ modulo $\lambda$ and is as well a prime integer, but no factor of $N (\eta-\eta_{r})$. Here, too, I \textit{no longer distinguish} the prime integers $p$'s from the $q$'s, but take the prime integers that belong to the exponent 1 just as falling \textit{under the same law}\footnote{The term ``Recht (right)'' is difficult to translate. To highlight the juridical connotation, I translate it as ``law''. However, I warn the reader that the translation might be misleading. Kummer simply pointed out the fact that they are all subject to the same jurisdiction.} as the others, the periods for these integers being simply with a single term and $f = 1$. (Die Bezeichnungen, welche ich habe, werden Ihnen wohl aus dem vorigen Briefe und sonst schon bekannt sein, n\"amlich $\lambda - 1 = e f$ und $q$ zum Exponenten $f$ geh\"orig modulo $\lambda$ und Primzahl aber kein Factor von $N (\eta-\eta_{r})$. Ich unterscheide auch hier gar nicht mehr Primzahlen~$p$ von den~$q$, sondern nehme die zum Exponenten 1 geh\"origen Primzahlen als eben \textit{in demselben Rechte} befindlich als die andern, es sind f\"ur diese nur die Perioden eingliedrig und $f = 1$.) \cite[p.~94, my emphasis]{46}.
\end{quote}

Just as in the first public presentation of 1846, Kummer's letter refers to the prime integers that are congruent to 1 modulo $\lambda$ using the letter $p$. Interestingly, in the letter, he refers to the other prime numbers using the letter $q$. Even though $p$ and $q$ are still used in the published article, $p$ still designates the first case considered, whereas $q$ now designates any prime integer (except $\lambda$), following the new insight that motivated the letter. This remark highlights the work that Kummer carried out to increase generality and uniformity in his treatment of divisibility, looking, as he later said, for a permanent property that would hold for all prime complex divisors and would thus allow him to introduce the ideal divisors in exactly the same fashion. The remark also sheds light on Kummer's discovery process. On this basis, let us thus return to how the 1846 article subsequently deals with the second definition for all prime complex divisors of prime integers. Kummer first makes a general statement about the periods:

\begin{quote}
One knows that the \textit{periods} are the $e$ roots of an \textit{equation} of degree $e$; and that this equation, interpreted as congruence, for the \textit{modulus $q$, always} has \textit{$e$ real roots of the congruence}, which I denote as $u, u_{1}, u_{2}, \dots, u_{e- 1}$, and which I take in a sequence corresponding to the periods ({\dots}). If, now, we denote simply with $\Phi(\eta)$ the \textit{complex number made up of periods} $c'\eta+ c'_1\eta_1+ c'_2\eta_2+ \dots + c'_{e-1}\eta_{e -1}$, \textit{then there are, among the prime numbers q} that belong to the exponent $f$, \textit{always some that} can be brought into the form $q = \Phi(\eta)\Phi(\eta_1)\Phi(\eta_2)\cdots \Phi(\eta_{e-1})$, in which the $e$ factors also never allow a further decomposition. (Bekanntlich sind die \textit{Perioden} die $e$ Wurzeln einer \textit{Gleichung} vom $e$ten Grade, und diese als Congruenz aufgefasst,\footnote{Here, instead of ``und diese als Congruenz aufgefasst,\dots'' \cite[p.~321]{25} has ``und diese, als Congruenz betrachtet, \dots (and this equation, considered as congruence, \dots)''.} f\"ur den Modul $q$, hat \emph{immer $e$ reale Congruenzwurzeln}, welche ich durch $u, u_1, u_{2}, \dots, u_{e- 1}$ bezeichne, und welche ich in einer entsprechenden Reihenfolge nehme wie die Perioden,\footnote{Here, \cite[p.~321]{25} reads as follows: ``und in einer entsprechenden Reihenfolge nehme, wie die Periode (and take in a corresponding sequence, like the periods)''.} ({\dots}). Wird nun \emph{die aus Perioden gebildete complexe Zahl} $c'\eta+ c'_1\eta_1+ c'_2\eta_2+ \dots + c'_{e-1}\eta_{e -1}$ kurz durch $\Phi(\eta)$ bezeichnet, so giebt \emph{es unter den Primzahlen} $q$, welche zum Exponenten $f$ geh\"oren, \emph{immer solche}, die sich auf die Form $q = \Phi(\eta)\Phi(\eta_1)\Phi(\eta_2)\cdots \Phi(\eta_{e-1})$ setzen lassen,\footnote{Here, \cite[p.~321]{25} has ``bringen lassen''.} in welcher auch die $e$ Factoren niemals eine weitere Zerlegung gestatten.) \cite[p.~90, my emphasis]{23}.
\end{quote}

In other words, Kummer established a link between the $e$ periods of length $f$ and $e$ integers with respect to a given prime integer $q$ that corresponds to the exponent $f$. He then turned to prime integers $q$ of this kind, whose decomposition into prime factors is actual (I use the term ``actual'' as a translation of \textit{wirklich} -- the term Kummer used above, when he referred to a decomposition into prime factors that can be carried out explicitly). The key point now for him, as it was in the first case, is to find a property of these divisors that does not depend on the ``actual'' decomposability, which he does in the subsequent sentences as follows:

\begin{quote}
Let us put, \textit{instead of the periods}, the \textit{roots} of the congruence that correspond to them ({\dots}). Hence one of the $e$ prime factors will always be congruent to~0 for the modulus~$q$. If, now, \textit{any complex number} $f(\alpha)$ contains the prime factor $\Phi(\eta)$,\footnote{In relation to the notation used above by Kummer for the decomposition of $q$, we should here rather have $\Phi(\eta_k$), to point out the fact that the divisor considered is the one corresponding to a circular permutation of the integers $u, u_1, u_2, \dots,u_{e-1}$, which is the order corresponding to $\eta,\eta_1,\eta_2,\dots,\eta_{e-1}$.\label{footnote89}} then it will have the \textit{property} that it becomes congruent to~0 for the modulus $q$, for $\eta=u_k$, $\eta_1=u_{k+1}$, $\eta_2=u_{k+2}$, \dots, etc. This property, now, [{\dots}] is a \textit{permanent property/property that holds also for the prime numbers} $q$ that \textit{do not allow a decomposition} into $e$ actual complex prime factors. This \textit{property could} \textit{hence} be used as \textit{definition} of the complex prime factors. \textit{However,} it would have the \textit{deficiency that} it would not express the equal ideal prime factors that occur in a complex number. (Setzt man \textit{anstatt}\footnote{Here, \cite[p.~321]{25} has ``statt''.} der \textit{Perioden} ihre entsprechenden \textit{Congruenzwurzeln}, ({\dots}) so wird immer einer der $e$ Primfactoren congruent Null, f\"ur den Modul $q$. Enth\"alt nun \textit{irgend eine complexe Zahl} $f(\alpha)$ den Primfactor $\Phi(\eta)$,\footnote{See footnote~\ref{footnote89}.} so wird sie die \textit{Eigenschaft} haben, f\"ur $\eta=u_k$, $\eta_1=u_{k+1}$, $\eta_2=u_{k+2}$, \dots, etc. congruent Null zu werden, f\"ur den Modul~$q$. Diese \textit{Eigenschaft} nun ({\dots}) ist eine \textit{bleibende} auch f\"ur \textit{diejenigen Primzahlen} $q$, welche \textit{eine Zerlegung} in die $e$ wirklichen complexen Primfactoren \textit{nicht gestatten}, sie \textit{k\"onnte daher als Definition} der complexen Primfactoren benutzt werden\footnote{In \cite[p.~321]{25}, Kummer modified the punctuation as follows: ``Diese $ $Eigenschaft$ $ nun (\dots ) ist eine bleibende: auch f\"ur diejenigen Primzahlen~$q$, welche eine Zerlegung in die $ $e$ $ wirklichen complexen Primfactoren nicht gestatten. Sie k\"onnte daher als Definition \dots''.} w\"urde aber \textit{auch} den \textit{Mangel haben}, dass sie die in einer complexen Zahl vorhandenen gleichen idealen Primfactoren nicht ausdr\"uckte.) \cite[p.~90, my emphasis]{23}.
\end{quote}

For any prime integer $q$ associated with the exponent $f$, for which a decomposition holds, Kummer highlighted a ``property (Eigenschaft)'' that is again stated as ``permanent'' and thus taken as ``defining'' the divisors whether the decomposition be actual or not. We thus see again the same conceptual tools as those introduced by Chasles, which reflects the same assumption: the ``ideal divisors'' are integral parts of the numbers they divide. What Kummer explained in the part of his October 18, 1845, letter quoted above highlights in which respect this way of dealing with any prime integer $q$ is simply a generalization of the first considerations about the~$p$'s.\looseness=-1

As the end of the quotation makes clear, once again this definition is not completely general, since it cannot encompass the prime divisors that might occur with a multiplicity strictly greater than 1 in the decomposition of a complex number into prime factors. For this, once again, Kummer looked for a still more general property that will hold for all cases. As he wrote:

\begin{quote}
The definition of ideal complex prime factors that I have chosen, which indeed essentially \textit{agrees} with the one that was indicated here, but is \textit{simpler} and \textit{more general}, rests on the fact, that ({\dots}here Kummer then introduces the permanent property that he will in the end adopt). It would lead me too far here, if I wanted to develop the relationship and the agreement between this definition and those that were indicated above, which were given through congruence conditions;{\dots} (Die von mir gew\"ahlte Definition der idealen complexen Primfactoren, welche im wesentlichen\footnote{Kummer \cite[p.~321]{25} capitalized the word.} zwar mit der hier angedeuteten \textit{\"ubereinstimmt}, aber zugleich\footnote{Kummer \cite[p.~321]{25} deleted this adverb. } \textit{einfacher} und \textit{allgemeiner} ist, beruht darauf, dass ({\dots}) Es w\"urde mich\footnote{Kummer \cite[p.~322]{25} deleted this pronoun.} hier zu weit f\"uhren, wenn ich den Zusammenhang und die \"{U}bereinstimmung dieser Definition mit den oben angedeuteten, welche durch Congruenzbedingungen gegeben werden, entwickeln wollte; {\dots}) \cite[p.~90--91, my emphasis]{23}.
\end{quote}

 We see that, here again, Kummer underlined the values of simplicity and generality on which he had laid stress in his introduction. I will not examine further the third definition that Kummer put forwards,\footnote{For the general idea, see \cite[p.~134--137]{18} and the previous proofs in which the key function $\Psi$ is introduced. However, here, Edwards did not strictly follow Kummer's procedure (see \cite[p.~129, footnote]{18}).} since I think the general structure of Kummer's procedure is clear: as he made explicit at the beginning of his development, he considered several permanent properties in his search for a property that would be uniform for all cases. The organization of the article puts forward in turn three such properties with an ever increasing generality, and at each of the first two steps, Kummer discarded the property considered on account of its lack of generality. However, at each step, the permanent property considered fulfilled the purpose that Chasles had outlined for properties of this kind: the property allowed Kummer to define at the same time actual and ideal divisors -- all then considered as integral components of the numbers they divide. We thus see clearly in which respect Kummer adopted the strategy suggested by Chasles to approach in a uniform way the type of elements that Poncelet had juxtaposed: those, like secants and chords, that had their mathematical identity ``actually'' and those who had it only ``ideally''.\looseness=-1

In his 1851 presentation of his theory in French, Kummer devoted a specific section to the ``Definition and general properties of ideal factors of a complex number (d\'efinition et propri\'et\'es g\'en\'erales des facteurs id\'eaux d'un nombre complexe)'' \cite[p.~425]{26}. After having defined the ``ideal factors'' of a prime integer $q$, Kummer commented on the value he attached to the introduction of factors of this kind, and he further explained how they compare with previous mathematical achievements. In this context, he returned for the second time to the parallel between his introduction of ideal factors and the fact of considering in geometry a straight line as the secant common to two circles whether the circles intersect or not. In correlation with the fact that he returned to this comparison, just as in 1846, the point of the comparison concerns the definition of ideal factors, and it aimed at justifying the definition chosen. Indeed, Kummer drew a parallel between his definition of ``ideal divisors'' and the fact of defining the common secant using the ``permanent'' property that characterizes the radical axis of two circles. Similarly, he equated the property that a complex number has an ``extant prime factor'' and the ``contingent'' property of the same line according to which it is the secant going through the points of intersection of the two circles \cite[p.~430]{26}. Despite the fact that the analogy with geometry is only one among three comparisons he made, only in this case does Kummer's comment bear on the justification of his definition on account of the fact that it, too, relies on a permanent property. This is the part of Kummer's work on complex numbers that is more specifically related to reflections in geometry developed by Chasles. A few lines earlier, Kummer had justified the ``denomination'' of ``ideal factor'' \cite[p.~429]{26}. In this respect, Kummer emphasized the deep analogy existing between the ordinary factors of a complex number and its ideal factors, just as Poncelet justified using the denomination of ``ideal secant'' to broaden the conceptions of geometry. Note that in both cases, the new terminology consists in adding an epithet -- and more importantly the same epithet: ``ideal'' -- to an ordinary concept (secant in geometry and divisor in the theory of divisibility) to establish a connection between cases that otherwise would remain perceived as different. Indeed, with the term ``ideal'', Poncelet designated a straight line whose connection with other parts of the figure have become imaginary, but which can nevertheless be constructed. By contrast, Kummer drew on Chasles's understanding when he asserted that ``the ideal factors make visible, so to say, the inner constitution of numbers (les facteurs id\'eaux rendent visibles, pour ainsi dire, la constitution int\'erieure des nombres)'' \cite[p.~429--430]{26}.\footnote{Kummer \cite[p.~92]{23} made a similar assertion.} In other words, the permanent properties shed light on what the integral components of a number are. These pieces of evidence from his 1851 paper confirm the part played by the ``new geometry'' in Kummer's reflections.\looseness=-1

\section{Conclusion}\label{section6}

This article has argued that Kummer's introduction of ``ideal divisors'' owed a debt to both Poncelet's concept of ``ideal elements'' and to Chasles's reconceptualization of the ``principle of continuity'' into the ``principle of contingent relationships''. Indeed, the use of the term ``ideal'' as well as the focus on ``permanent properties'' to ``define'' ``ideal divisors'' are central facets of the 1846 public presentation of Kummer's ``theory of complex numbers'', and they allow us to interpret the meaning of Kummer's explicit parallel between his work and recent work in the ``new geometry''.

In addition to establishing a historical fact, what does this conclusion tell us about Kummer's ideal divisors?

I have emphasized above that reading Kummer's first public presentation of his theory of complex divisors against the backdrop of Chasles's \textit{Aper\c{c}u historique} highlights the part played in Kummer's reflections by Chasles's concepts of ``contingent'' and ``permanent properties'', as well as by the French geometer's urge to look for definitions that would rest on permanent properties. However, the very fact that Kummer embraced these ideas sheds light on his conception of these ``ideal divisors''. For him, like the ideal secant common to the two circles, ``ideal divisors'' are ``integral'' components of the numbers they divide: this is what justifies looking for a permanent property that will define them in all circumstances. This conclusion fits with the claim Kummer made at the beginning of his 1846 presentation, when he asserted that the ``prime factors'' of complex numbers are ``true (usually ideal)''.

The only difference between Kummer and Chasles in this respect is that Kummer keeps referring to the permanent factors that he seeks to identify -- as well as to the straight line that persists when the circles no longer cut each other -- using the term ``ideal''. This, as well as the fact that Kummer makes clear that for him these ``prime factors'' are ``true (usually ideal)'', suggests a similarity between Poncelet's and Kummer's understanding of these permanent elements. Indeed, for both of them these elements are different from ``imaginary'' elements. We recall that, for Poncelet, ``imaginary'' elements are those that cannot be constructed. In contrast, for him, ``ideal'' elements are parts of a figure that would subsist in any general position of that figure. Moreover, although the relationships between the ``ideal elements'' and other parts of the figure may have become ``imaginary'', one can find ways to construct these ``ideal elements''. In this respect, it is interesting to note that Kummer's ideal numbers are also only defined not in and of themselves, but only \textit{by their relations} with other complex numbers. Moreover, it is quite remarkable that what Kummer did is precisely to provide a way to construct these ideal divisors using congruences.

Here, we may go one step further. Indeed, as early as in the October 18, 1845 letter, Kummer ``justified'' the use of the expression ``ideal prime factor'' by reference to the fact that these factors make ``computation'' with complex numbers identical to computations with integers. This reference to what we may call, using Dedekind's later expression, ``the general laws of divisibility''\footnote{See, e.g., \cite[p.~278]{14}: ``les lois g\'en\'erales de la divisibilit\'e''.} evoke the full name that Poncelet had given to his key principle: ``the \textit{principle} of permanence or \textit{indefinite continuity of the mathematical laws of magnitudes.}'' In both cases, what is at stake is a persistence of the mathematical relationships between the entities~-- a~persistence that takes on the guise of total uniformity. This fact is correlated with the emphasis Kummer places on the values of simplicity and generality, which were also crucial epistemological values for Chasles~\cite{9}.

The preceding remarks explain in which respect we might consider that Kummer realized a~synthesis between Poncelet's and Chasles's approach to ``ideal'' elements. More precisely, Kummer's achievement drew on a synthesis of the philosophical reflections underpinning the shaping of the new geometry in the hands of geometers such as Poncelet and Chasles. Indeed, their philosophical reflections were instrumental in fashioning what would later become projective geometry \cite{9}. Moreover, it was precisely these reflections that were picked up from geometry into the ``theory of complex numbers'', contributing to a major change in this latter context.

\subsection*{Acknowledgements}

It is a pleasure to dedicate this article to Jean-Pierre Bourguignon, without whom I might never have known that the mathematics section of CNRS intended to open hirings to the history of mathematics. Without this piece of information, I might never have continued research in this domain. Given this context, I find it appropriate to dedicate to him an article showing how the historical and philosophical reflections developed by mathematicians led to significant advancements in mathematics. I have been interested in the history of geometry in France in the first half of the 19${}^{\rm th}$ century since the 1980s, thanks to Serge Pahaut's (1945--2023) invitation to work with him on duality. His passing away just as I was completing this article infuses great sadness with some of its pages. My interest in the notions of ideality was spurred by the ANR project ``Ideals of proof'' that Mic Detlefsen (1948--2019) developed at the University Paris Diderot between 2007 and 2011. In this context, I gave a joint talk with Bruno Belhoste on Poncelet's ideal elements (titled: ``Poncelet's ideal elements in geometry: between Carnot and Chasles'', March~5, 2009), which we intend to publish. The inspiration for the present article derives from my first encounter with Ernst Kummer's writings on number theory, in the context of a talk given at our research group (REHSEIS at the time, SPHERE, now) by Jacqueline Boniface on the topic. As we incited speakers to do at the time, she handed out pages of Kummer's original text. I immediately recognized that these pages had been inspired by Michel Chasles's ``principle of contingent relationships''. This formed part of the results that we published in~\cite{11}, which concluded the research project carried out between 2004 and 2016, about the epistemic and epistemological value of generality (about these expressions, I refer the reader to the introduction of the book). In the book~\cite{9} I analyzed Chasles's understanding of generality in geometry, and mentioned the impact of his reflections, notably those relating to the aforementioned principle, on Kummer's work in number theory. In~\cite{4}, Boniface described Kummer's creation of ideal numbers, mentioning the analogy with geometry. The prologue to the book~\cite{10} took this case as an example of a circulation of scientists' reflections about an epistemological value from one context to another context. The present article builds on all these previous researches, and I am in particular grateful to Jacqueline Boniface for having introduced me to Kummer's work in number theory. I am also immensely grateful to Harold Edwards, since, without his book~\cite{18}, I might not have been able to enter into the maze of Kummer's mathematical world. Last, but not least, since 2022, Bruno Belhoste and I have begun a joint project on Poncelet's works. Our common work has helped me refine my understanding of Poncelet's philosophical and geometrical reflections. Since I completed the first version of this article, I was fortunate to receive remarks from Christian Houzel, Max Lindh, Colin McLarty, Patrick Popescu-Pampu, Jean-Pierre Serre as well as from three referees. I thank them all for their comments which have enabled me to significantly improve the first version of the article. The contribution of Jeremy Gray and Oussama Hijazi to the final version of this article was of crucial importance.


\pdfbookmark[1]{References}{ref}
\LastPageEnding

\end{document}